%

%


\def\today{\ifcase\month\or January\or February\or
March\or April\or May\or June\or July\or August\or
September\or October\or November\or December\fi
\space\number\day, \number\year}




\def\dspace{\lineskip=2pt\baselineskip=18pt
\lineskiplimit=0pt}

\font \bbrm=cmbx10 at 12pt

\def\bigtype{\bbrm}

\hsize=13.5cm
\magnification=1200
\def\ce{\centerline}

\def\hb{\hfill\break}

\def\title #1{\null\bigskip\ce{\bigtype #1}
\bigskip}

\def\alp{\alpha}		
\def\bet{\beta}		
\def\gam{\gamma}		
\def\del{\delta}		\def\Del{\Delta}
\def\eps{\varepsilon}		

		\def\Tet{\Theta}

\def\kap{\kappa}
\def\lam{\lambda}		
\def\sig{\sigma}		

\def\ome{\omega}		


\def\calD{{\cal D}}

\def\calF{{\cal F}}

\def\calK{{\cal K}}

\def\calP{{\cal P}}



    
\font\tenboldgreek=cmmib10
 \font\sevenboldgreek=cmmib10 at 7pt
\font\fiveboldgreek=cmmib10 at 7pt
\newfam\bgfam
\textfont\bgfam=\tenboldgreek
\scriptfont\bgfam=\sevenboldgreek
\scriptscriptfont\bgfam=\fiveboldgreek

\mathchardef\ggarrow="7010

\font\tengerman=eufm10 \font\sevengerman=eufm7
\font\fivegerman=eufm5
\font\tendouble=msym10 \font\sevendouble=msym7
\font\fivedouble=msym5

\textfont4=\tengerman \scriptfont4=\sevengerman
\scriptscriptfont4=\fivegerman
\newfam\dbfam
\textfont\dbfam=\tendouble \scriptfont\dbfam=
\sevendouble
\scriptscriptfont\dbfam=\fivedouble

\mathchardef\ng="702D
\mathchardef\dbA="7041
\mathchardef\sm="7072
\mathchardef\nvdash="7030
\mathchardef\nldash="7031
\mathchardef\lne="7008
\mathchardef\sneq="7024
\mathchardef\spneq="7025
\mathchardef\sne="7028
\mathchardef\spne="7029
\mathchardef\ltms="706E
\mathchardef\tmsl="706F

\mathchardef\dbA="7041


\mathchardef\dbA="7041 
\mathchardef\dbB="7042 
\mathchardef\dbC="7043 
\mathchardef\dbD="7044 
\mathchardef\dbE="7045 
\mathchardef\dbF="7046 
\mathchardef\dbG="7047 
\mathchardef\dbH="7048 
\mathchardef\dbI="7049 
\mathchardef\dbJ="704A 
\mathchardef\dbK="704B 
\mathchardef\dbL="704C 
\mathchardef\dbM="704D 
\mathchardef\dbN="704E 
\mathchardef\dbO="704F 
\mathchardef\dbP="7050 
\mathchardef\dbQ="7051 
\mathchardef\dbR="7052 
\mathchardef\dbS="7053 
\mathchardef\dbT="7054 
\mathchardef\dbU="7055 
\mathchardef\dbV="7056 
\mathchardef\dbW="7057 
\mathchardef\dbX="7058 
\mathchardef\dbY="7059 
\mathchardef\dbZ="705A 

\def\nek{,\ldots,}
\def\sdp{\times \hskip -0.3em {\raise 0.3ex
\hbox{$\scriptscriptstyle |$}}} 


\def\dom{\mathop{\rm dom}\nolimits}

\def\min{\mathop{\rm min}}
\def\MOD{\mathop{\rm mod}}



\def\oF{{\overline F}}

\def\oG{{\overline G}}

\def\op{{\overline p}}

\def\oq{{\overline q}}

\def\OR{{\overline r}}

\def\oV{{\overline V}}


\def\otau{{\overline\tau}}







\def\ddownarrow{\big\downarrow \hskip-0.70em\raise
2pt\hbox {$\big\downarrow$}}
\def\longright #1#2 {\smash{\mathop{\hbox to
#1pt {\rightarrowfill}}\limits_{#2}}}
\def\sqr#1#2{{\vcenter{\hrule height.#2pt\hbox{\vrule
width.#2pt height#1pt \kern#1pt \vrule width.#2pt}
\hrule height.#2pt}}}
\def\square{\mathchoice{\sqr34}{\sqr34}{\sqr{2.1}3}
{\sqr{1.5}3}}

\def\buildrul#1\under#2{\mathrel{\mathop{\null#2}
\limits_{#1}}}

\def\boxit#1{\vbox{\hrule\hbox{\vrule\kern3pt
\vbox{\kern3pt#1 \kern3pt}\kern3pt\vrule}\hrule}}

\def\subheading#1{\medskip\goodbreak\noindent{\bf
#1.}\quad}

\def\sect#1{\goodbreak\bigskip\centerline{\bf#1}
\medskip}
\def\pr{\smallskip\noindent{\bf Proof:\quad}}
\def\onumber #1{\ooalign{\hfil\raise.07ex\hbox{
\hfill$\scriptstyle \,#1$\hfil}
\cr\cr{$\bigcirc$}}}
\def\onumber c{\ooalign{\hfil\raise.07ex\hbox
{\hfill$\scriptstyle \,c$\hfil}
\cr\cr{$\bigcirc$}}}
\def\alpcirc {\ooalign{\hfil\raise.07ex
\hbox{\hfill$\scriptstyle\alp\;$\hfill}\cr\cr
{$\bigcirc$}}}

\def\longmapright #1#2 {\smash{\mathop{\hbox to
#1pt {\rightarrowfill}}\limits^{#2}}}
\def\longmapleft #1 #2 {\smash{\mathop{\hbox to
#1 pt {\leftarrowfill}}\limits^{#2}}}

\def\references#1{\goodbreak\bigskip\par\centerline
{\bf References}\medskip\parindent=#1pt}
\def\ref#1{\par\smallskip\hang\indent\llap{\hbox
to \parindent{#1\hfil\enspace}}\ignorespaces}

\def\back{{\raise 2.5pt\hbox{$\,\scriptscriptstyle
\backslash\,$}}}
\def\bks{{\backslash}}
\def\part{\partial}
\def\lwr #1{\lower 5pt\hbox{$#1$}\hskip -3pt}
\def\rse #1{\hskip -3pt\raise 5pt\hbox{$#1$}}
\def\lwrs #1{\lower 4pt\hbox{$\scriptstyle #1$}
\hskip -2pt}
\def\rses #1{\hskip -2pt\raise 3pt\hbox
{$\scriptstyle #1$}}

\def\<#1{\left\langle{#1}\right\rangle}

\def\subinbn{{\subset\hskip-8pt\raise 0.95pt
\hbox{$\scriptscriptstyle\subset$}}}

\def\llvdash{\mathop{\|\hskip-2pt
\raise 3pt\hbox{\vrule height 0.25pt width 1.5cm}}}

\def\lvdash{\mathop{|\hskip-2pt \raise 3pt\hbox
{\vrule height 0.25pt width 1.5cm}}}

\def\fakebold#1{\leavevmode\setbox0=\hbox{#1}%
  \kern-.025em\copy0 \kern-\wd0
  \kern .025em\copy0 \kern-\wd0
  \kern-.025em\raise.0333em\box0 }

\font\msxmten=msxm10
\font\msxmseven=msxm7
\font\msxmfive=msxm5
\newfam\myfam
\textfont\myfam=\msxmten
\scriptfont\myfam=\msxmseven
\scriptscriptfont\myfam=\msxmfive
\mathchardef\rhookupone="7016
\mathchardef\ldh="700D
\mathchardef\leg="7053
\mathchardef\ANG="705E
\mathchardef\lcu="7070
\mathchardef\rcu="7071
\mathchardef\leseq="7035
\mathchardef\qeeg="703D
\mathchardef\qeel="7036
\mathchardef\blackbox="7004
\mathchardef\bbx="7003
\mathchardef\simsucc="7025

\def\rhookup{{\fam=\myfam \rhookupone}}

\def\bigsquare{{\fam=\myfam\bbx}}

\font\tencaps=cmcsc10
\def\smallcaps{\tencaps}

\def\author#1{\bigskip\ce{\smallcaps #1}\medskip}

\def\upddots{\mathinner{\mkern
1mu\raise 1pt \hbox{.}\mkern 2mu \mkern
2mu \raise 4pt\hbox{.}\mkern 1mu \raise 7pt\vbox
{\kern 7 pt\hbox{.}}} }

\def\varchi{\ooalign{{\raise
1.385pt\hbox{$\chi$}}\crcr\hbox{--}\crcr}}

\def\trianarrow{{\raise 2pt\hbox to 0.50cm
{\hrulefill}\triangleright}}

\null
\overfullrule=0pt
\def\vline{\Vert\hskip -1.4pt
{\vcenter{\hrule width 1truecm}}}  
\def\ve{|\!\!\!=}
\def\su{\buildrul \sim\under}
\def\aU{\vec U}
\def\aV{\vec V}
\def\aC{\vec C}
\def\aF{\vec\calF}
\def\Col{\rm Col}
\sect{CARDINAL PRESERVING IDEALS}

\vskip 1truecm
\ce{by}
\vskip 1truecm
$$\vbox{\halign{\tabskip 2em \hfil#\hfil&\hfil#\hfil&\hfil
#\hfil\cr
Moti Gitik &&Saharon Shelah\cr
School of Mathematical Sciences&and&Department
of Mathematics\cr
Tel Aviv University&&Hebrew University of
Jerusalem\cr
Tel Aviv, Israel&& Jerusalem, Israel\cr}}$$
\vskip 2truecm
\dspace
\sect{Abstract}

\noindent
We give some general criteria, when
$\kap$-complete forcing preserves largeness
properties -- like $\kap$-presaturation of
normal ideals on $\lam$ (even when they
concentrate on small cofinalities).  Then we
 quite accurately obtain the consistency
strength ``$NS_\lam$ is
$\aleph_1$-preserving", for $\lam >
\aleph_2$.
\vfill\eject

We consider the notion of a $\kap$-presaturated
ideal which was basically introduced by
Baumgartner and Taylor [B-T].  It is a
weakening of presaturation.  It turns out
that this notion can be preserved under
forcing like the Levy collapse.  So in
order to obtain such an ideal over a small
cardinal it is enough to construct it over
an inaccessible and then just to use the
Levy collapse.

The paper is organized as follows.  In
Section 1 the notions are introduced and
various conditions on forcing notions for
the
preservation of $\kap$-presaturation are
presented.  Models with $NS_\lam$ cardinal
preserving are constructed in Section 2.
The reading of this section requires some
knowledge of [G1,2,4].

The results of the first section are due to
the second author and second section to the
first.

\subheading{Notation}  $NS_\kap$ denotes
the nonstationary ideal over a regular
cardinal $\kap > \aleph_0$, $N S_\kap^\mu$
denotes the $NS_\kap$ restricted to
cofinality $\mu$, i.e. $\{X\subseteq
\kap|X\cap\{\alp <\kap| {\rm cf}\
\alp=\mu\} \in NS_\kap\}$, $D$ denotes a normal
filter over a regular cardinal $\lam =
\lam(D) > \aleph_0$, $D^+ = \{A \subseteq
\lam|A \ne \emptyset$ mod $D$, i.e. $\lam -
A \not\in D\}$. $\buildrul \sim\under D^Q$
for a forcing notion $Q$, such that 
$\buildrul Q \under \vline$ ``$\lam(D)$ is regular" is the
$Q$-name of the normal filter 
on $\lam(D)$ generated by $D$. By forcing
with $D^+$ we mean the forcing with
$D$-positive sets ordered by inclusion.

\sect{1.~~$\kap$-Presaturation and Preservation
Conditions}

The definitions and the facts 1.1-1.5 below
are basically due to Baumgartner-Taylor
[B-T].
\subheading{Definition 1.1}  A normal
filter (or ideal) $D$ over $\lam$ is
$\kap$-presaturated if
$\buildrul {D^+} \under 
\vline$  
``every set of ordinals of cardinality $<
\kap$ can be covered by a set of $V$ of
cardinality $\lam^{''}$.

Let us formulate an equivalent definition
which is much  easier to use.

\subheading{Definition 1.2}  A normal
filter $D$ on $\lam$ is
$\kap$-presaturated, if for every $\kap_0<
\kap$, and maximal antichains $I_\alp =
\{A_i^\alp: i < i_\alp\}$, i.e.
$$A_i^\alp \in D^+, [i \ne j \Rightarrow
A_i^\alp \cap A_j^\alp \not\in D^+], \quad
(\forall A \in D^+)(\exists i < i_\alp) A
\cap A_i^\alp \in D^+)$$
(for $\alp < \kap_0$) for every $B \in D^+$
there is $A^* \in D^+$ such that $A^*
\subseteq B$ and $(\forall \alp < \kap_0)|
\{ i < i_\alp: A^* \cap A_i^\alp \in D^+\}|
\le \lam$.

\proclaim Fact 1.3. If $\kap < \lam$,
$(\forall \Tet < \lam) [\Tet^\kap < \lam]$,
$\{i < \lam: cfi\ge \kap\} \in D$, $D$ is
$\kap$-presaturated, then in 1.2 we can
find $A^{**}$ such that $\forall \alp <
\kap \exists ! i < i_\alp$ \enskip $A^{**}
\cap A_i^\alp \in D^+$.

\pr First pick $A^*$ as in the definition.
For every $\alp < \kap_0$ let $\langle
B_r^\alp| \tau < \lam_\alp \le \lam\rangle$
be an enumeration of $\{A_i^\alp \cap A^*|
i < i_\alp, \quad A_i^\alp \cap A^* \in
D^+\}$.  Without loss of generality min
$B_r^\alp > \tau$.  For $\nu< \lam$, $\alp
< \kap_0$ let $\tau_\alp(\nu)$ be the 
least $\tau$ s.t. $\nu \in B_\tau^\alp$ if
such $\tau$ exists and -1 otherwise.
Define $f:\lam \to \lam$ by $f(\nu) =
\bigcup\limits_{\alp < \kap_0}
\tau_\alp(\nu)$.  Then there are $\Tet <
\lam$ and a $D$-positive subset $A'$ of
$A^*$ s.t $f^{''} (A') = \{\Tet\}$.  Since
$\Tet^\kap < \lam$, using once more the
Fodor Lemma, we can find $A^{**} \subseteq
A'$ as required.

\subheading{Definition 1.4} A normal ideal
or filter on $\lam$ is $\kap$-preserving
iff
\item{(a)} $I$ is precipitous
\item{(b)} $\buildrul I^+ \under \vline$ ``$\check \kap$ is a
cardinal".
 
\subheading{Remarks}

\item{(1)} If $\kap = \lam^+$, then such
an ideal is called presaturated.  This notion
was introduced by Baumgartner-Taylor [B-T].
\item{(2)} It is unknown if for $\kap \ge
\ome_1$ (b) may hold without (a).  But if
$2^\lam = \lam^+$, then (b) $\Rightarrow$
(a).
\item{(3)} Every precipitous ideal is
$|I^+|^+$-preserving.

\proclaim Proposition 1.5.  Suppose that
$\kap<\lam$ are regular, $2^\lam = \lam^+$
and $D$ is a normal ideal over $\lam$.
Then the following conditions are
equivalent:
\item{(a)} $D$ is not $\kap$-presaturated;
\item{(b)} $\buildrul D \under \vline$ cf $\lam^+ < \kap$;
\item{(c)} the forcing with $D^+$ collapses
all the cardinals between $\lam^+$ and some
cardinal $<\kap$;
\item{(d)} $I$ is not $\kap$-preserving;
\item{(e)} a generic ultrapower of $D$ is
not well founded or it is well founded but
it is not closed in $V^{D^+}$ under less
than $\kap$-sequences of its elements;\hb
If in addition $2^{<\kap} < \lam$ then also
\item{(f)} the forcing with $D^+$ adds new
subsets to some ordinal $< \kap$.

\pr (b) $\Rightarrow$ (c) since $|D^+| =
\lam^+$.  The direction (c) $\Rightarrow$
(a), (b), (d), (e), (f) are trivial.  If
$\neg$ (b) holds, then, as in [B-T], also
$\neg$(a) holds.  $\neg$(a) implies
$\neg$(d), $\neg$(e) and $\neg$(f).

We are now going to formulate various
preservation conditions.

\proclaim Lemma 1.6.  
Suppose that $\kap < \lam$ are regular
cardinals, $D$ a normal filter on $\lam =
\lam(D)$, $Q$ a forcing notion 
 $\buildrul Q \under \vline$ ``$\lam$ regular" and
$\su D^Q$ denotes the normal filter on $\lam$
in $V^Q$ which $D$ generates. \hb
If $D$ is $\kap${\it -presaturated} (in
$V$) then $\su D^Q$ is $\kap$-presaturated in $V^Q$
provided that for some $K$
\hb
$(*)^1_{D,K,Q,\kap}$ $K$ is a structure
with universe $|K|$, partial order $\le =
\le^K$ unitary function $T = T^K$, and
partial unitary functions $p_i = p_i^K$
such that
\item{(a)} for $t \in K$, $T(t) \in D^+$;
\item{(b)} $K\ve t \le s$ implies $T(s)
\subseteq T(t) \MOD D$;
\item{(c)} $p_i(t)$ is defined iff $i \in
T(t)$;
\item{(d)} $p_i(t) \in Q$ and $[K\ve t\le
s$, $i \in T(s) \Rightarrow p_i(t) \le p_i
(s) ]$;
\item{(e)} $\buildrul \sim\under\tau_t=
\buildrul \sim\under\tau_t^K = \{i \in
T(t): p_i(t) \in \buildrul\sim \under G_Q\}
\subseteq \lam$ is not forced to be
$\emptyset$ mod $\buildrul \sim\under D^Q$;
\item{(f)} if $q\buildrul Q \under \vline
^{''}\buildrul \sim\under T \in
(\buildrul\sim\under D^Q)^{+\ ''}$ for
some $q \in Q$ {\it then} there is $s$, $s
\in K$, $p_i(s) \ge q$ and $\buildrul Q\under\vline^{''}
\buildrul {\sim_s} \under \tau \subseteq\su
T \MOD \su D^{Q^{''}}$;
\item{(g)} if $\kap(*) < \kap$, for $\alp <
\kap(*)\enskip Y_\alp \subseteq K$ is such that
$\{T(t): t \in Y_\alp\}$ is a maximal
antichain of $D^+$ such that $[t \in
Y_\alp, \quad s \in Y_\bet, \quad \alp >
\beta \Rightarrow (\exists C \in D)
(\forall i \in C \cap T(t) \cap T(s))
p_i(s) \le p_i(t)]$, and $T \in D^+$ {\it
then\/} there is $s$, $s \in K$, $T(s)
\subseteq T$ and there are $y_\alp
\subseteq Y_\alp$, $|y_\alp| \le \lam$ for
$\alp < \kap(*)$ such that $[\alp <
\kap(*), \quad x \in Y_\alp\backslash
y_\alp \Rightarrow T(x) \cap T(s) =
\emptyset \MOD D]$ and $\forall \alp <
\kap(*) (\forall x \in y_\alp) (\exists C
\in D) (\forall i \in C \cap T(x) \cap
T(s)) [p_i(x) \le p_i(s)]$.

\pr Let $\kap(*) < \kap$, and $
\su I_\alp = \{\su A^\alp_i: i < \su
i_\alp\}$ 
(for $\alp < \kap(*)$) be $Q$-names of
maximal antichains of $(\su D^Q)^+$ (as in
Definition 1.2) (i.e. $\buildrul Q \under
\vline^{''} \su I_\alp$ is a maximal
antichain".

We define by induction on $\alp < \kap(*)
\enskip Y_\alp h_\alp$ such that
\item{(i)} $Y_\alp \subseteq K, \quad
\{T(t): t\in Y_\alp\}$ a maximal antichain
of $D^+$;
\item{(ii)} for every $\beta<\alp$, $t \in
Y_\beta$, $s \in Y_\alp$ for some $C \in D
\enskip i \in C \cap T(s) \cap T(t) \to p_i(t) \le
p_i(s))$;
\item{(iii)} $h_\alp:Y_\alp \to$ ordinals,
and $\vline^{''} \tau_t \subseteq \su
A_{h_\alp (t)}^{\alp^{''}}$ for $t \in
Y_\alp$.

Using (g) pick $s$ and $\langle y_\alp|\alp
< \kap(*) \rangle$.  Then $\su \tau_s$ will
be as required in Definition 1.2.

\subheading{Definition 1.7}

\item{1)} 
$Gm(D,\gam,a)$ (where $a \subseteq \gam)$
is a game which lasts $\gam$ moves, in the
$i$th move if $i \in a$ player I, and if $i
\not\in a$ player II, choose a set $A_i \in
D^+$, $A_i \subseteq A_j \MOD D$ for $j <
i$.  If at some stage there is no legal
move, player I wins, otherwise player II
wins.

\item{2)} If we omit $a$, it means $a = \{1
+ 2i: 1 + 2i< \gam\}$.
\item{3)} $Gm^+ (D,\gam, a)$ is defined
similarly, but for player II to win, 
$\cap A_i \ne \emptyset$ has to hold as
well.

\medskip
By Galvin-Jech-Magidor [G-J-M], the
following holds.

\proclaim Proposition 1.8. $D$ is
precipitous if player I has no winning
strategy in $Gm^+ (D,\omega)$.

\proclaim Lemma 1.9.
\item{1)}  Suppose that $\gam<\lam$ are
regular cardinals, $D$ a normal filter on
$\lam = \lam(D)$, $Q$ a forcing notion
$\buildrul Q \under \vline^{''}\lam$
regular$''$, $\lam = \lam(D)$.
\hb
If player II wins the game $Gm(D,\gam, a)$
in $V$, then it wins in $Gm(D,\gam,a)$ in
$V^Q$ provided that for some $\kap$, $\gam
\le \kap$, the following principle holds:\hb
$(*)^2_{D,K,Q}$:
\hb
$K$ is a structure as in Lemma 1.6,
\itemitem{(a)} for $t \in K$, $T(t) \in
D^+$;
\itemitem{(b)} $K \ve t \le s$ implies $T(s)
\subseteq T(t) \MOD D$;
\itemitem{(c)} $p_i(t)$ is defined iff $i \in
T(t)$;
\itemitem{(d)} $p_i(t) \in Q$ and $[K \ve t \le
s,\quad i \in T(s) \Rightarrow p_i(t) \le
p_i(s)]$;
\itemitem{(e)} $\su \tau_t = \{i \in T(t):
p_i(t) \in \su G_Q\}$ and $p_i(t) \buildrul
Q \under \vline^{''} \su \tau_t \in (\su
D^Q)^{+\ ''}$ for $i \in T(t)$;
\itemitem{(f)} if, $q\buildrul Q \under
\vline^{''} \su T \in (\su D^Q)^{+\ ''}$ for
some $q \in Q$ {\it then} there is $s$, $s
\in K$, $p_i(s) \ge q$ and $\buildrul Q
\under \vline^{''} \su \tau_s \subseteq \su
T \MOD \su D^{Q\ ''}$;
\itemitem{(g)} $K$ is $\kappa$-complete, i.e.
if $\langle t_i: i < \del < \kap\rangle$ is
increasing and $\cap T(t_i) \in D^+$ then
there is $t \in K$, $t_i \le t$;
\itemitem{(h)} $Q$ is $\kap$-complete and
preserves stationarity of all $A \in D^+$.
\item{2)} The same holds with Player I
having no winning strategy.
\vfill\eject

\pr
\item{1)}  
Let us provide a winning strategy for
player II.  Let $G \subseteq Q$ be generic
over $V$ without loss of generality the
players choose $Q$-names for their moves.
We will now describe the strategy of player
II in
$V[G]$.
\hb
Player II also chooses $T_i, t_j(j \in
a, i < \gam)$, according to the moves.
Player II preserves the following (for a
fixed winning strategy $\oF$ of player II
in $Gm(D, \gam, a)$ in $V$).
\item{(*)} the plays $\langle \su A_i: i <
\gam^* \le \gam\rangle$ and $\langle T_i: i
< \gam^*\rangle$ so far satisfy
\itemitem{(a)} $\langle T_i: i <
\gam^*\rangle$ is a beginning of a play of
$Gm(D,\gam,a)$ (in $V$) in which player II
uses the strategy $\oF$;
\itemitem{(b)} for $i \in a$, $\su A_i =
\su \tau_{t_i}$, $T(t_i) = T_i$, $t_i \in
K$;
\itemitem{(c)} for $j < i < \gam^*$, $t_j
\le t_i$;
\itemitem{(d)} $V[G]\ve ''\su A_i [G] \in
D^{+\ ''}$.

\proclaim Proposition 1.10.  Suppose $\kap
< \lam$ are regular cardinals, $\forall
\Tet < \lam[\Tet^{<\kap} < \lam],\ D$ a normal
filter on $\lam, Q$ a forcing notion,
$\buildrul Q \under \vline^{''} \lam$ is a
regular$''$ and $\su D^Q$ denotes the normal
filter on $\lam$ in $V^Q$ which $D$
generates.  Let\hb
$(*)^3_{D,K,Q}$ denote the following:
\itemitem{(a)} $K$ is a partial order,
elements of $K$ are of the form $\op =
\langle p_i: i \in T\rangle$, where $T =
T_{\op}$ is $D$-positive and $p_i \in Q$
for $i \in T$.  For $\op, \oq \in K$, $\op
\ge \oq$ iff $T_{\op} \subseteq T_\oq$ and
for every $i \in T_\op$, $p_i \ge q_i$ in
the ordering of $Q$;
\itemitem{(b)} if $\op \in K$, then
$\tau_\op = \{ i \in T: p_i \in \su G_Q\}$
is not forced to be empty mod $D$;
\itemitem{(c)} if $\su T$ is a $Q$-name, $\su T
\cap \tau_\op$ not forced to be empty mod
$D$, then for some $\oq \in K$, $\oq \ge
\op$ and $\buildrul Q \under \vline^{''}
\tau_\oq \subseteq \su T^{''}$;
\itemitem {(d)} (I) if $\op^\alp_i (\alp <
\alp(*) < \kap)$ is an increasing sequence
in $K$ and $\bigcap\limits_{\alp<\alp(*)}
T^\alp \in D^+$ then there is an upper
bound in $K$ or (II) if $\op^\alp (\alp <
\alp(*) < \kap)$ is a sequence from $K, T
\in D^+$, $T =
\bigcap\limits_{\alp<\alp(*)} T^\alp$,
$[\alp < \bet$, $i \in T^\alp \cap T^\beta
\Rightarrow p_i^\alp \le p_i^\beta]$ then
there is $\op^* \in K$, $T_{\op^*} \subseteq
T$, $[i \in T^\alp \cap T_{\op^*}
\Rightarrow p_i^\alp \le p_i^*]$;
\itemitem{(e)} for every $g \in Q$, $\su
\tau$ s.t. $q \vline^{''} \su \tau \in
(D^Q)^{+\ ''}$ there is $\op \in K$, $p_i \ge
q$ and $q \vline \tau_\op \subseteq \tau$;
\itemitem{(f)} $\langle\emptyset_Q| i <
\lam \rangle \in K$, where $\emptyset_Q$ is
the minimal element of $Q$.

{\sl If $(*)^3_{D,K,Q}$ holds then also the
following statements are true: 
\item{(1)} If $D$ is $\kap$-presaturated in
$V$, $\{\del: cf \del \ge \kap$ (in $V)\}
\in D$, {\it then} $\su D^Q$ is
$\kap$-presaturated in $V^Q$.
\item{(2)} If player II wins in $Gm(\lam,
D, \gam)$ in $V$ $(\gam \le \kap)$, $Q$
$\kap$-complete then he wins in $Gm(\lam,
D, \gam)$ in $V^Q$ provided that

$$(*)^4_{D,K,Q}: \matrix{
(a)&\op \in K \Rightarrow p_i \buildrul Q
\under \vline^{''} \su \tau_\op \ne
\emptyset \MOD \su D^{Q\  ''} \ {\rm for}\ i
\in T_\op\cr
(b)&(*)^3_{D,K,Q}\cr}$$
\item{(3)} The same holds if we replace
``The player II wins" by ``Player I does
not have winning strategy" or the game
$Gm(\lam, D, \gam)$ is replaced by $Gm^+
(\lam, D, \gam)$.
\item{(4)} In (2) we can replace
$(*)^3_{D,K,Q}$ by:
\hb
$D_1$ a normal filter over $\lam$ in $V,
K_1$ like $K$,\hb
player II wins also in $Gm(\lam, D, \gam)$
and

$$(*)^5_{D,K,D_1,K_1,Q}: \matrix{
(\alp)&(*)^4_{D_1,K_1,Q}\cr
(\bet)&(*)^4_{D,K,Q}\cr
(\gam)&\hbox{ for every}\ \op = \langle p_i:
i \in T_\op \rangle \in K\ \hbox{there
are}\cr
&\oq = \langle q_\alp: \alp \in S \rangle
\in K_1 \quad \op^\alp = \langle p_i^\alp:
i \in T_\alp\rangle \in K\cr
&T_\alp \subseteq T_\op\ \hbox{pairwise
disjoint}\cr
&[i \in T_\alp \Rightarrow p_i \le
p_i^\alp], \quad q_\alp \vline^{''} \su
\tau_{\op^\alp} \ne \emptyset \MOD \su
D^{Q\ ''}\ .\cr}$$
}

\pr Let us prove (1).  The proof of (2) is
similar.  Let $\kap(*) < \kap$, $\su I_\alp
= \{\su A_i^\alp: i< \su i_\alp\}$ $(\alp <
\kap(*))$ be $\kap(*)$ $Q$-names of maximal
antichains of $D^+$ as mentioned in the
definition of $\kap$-presaturitivity (i.e.
it is forced that they are like that).

We define by induction on $\alp\le
\kap(*)$, $Y_\alp$, $j(\nu)$ and $\op^\nu
(\nu \in Y_\alp)$ such that
\item{(i)} $Y_\alp$ is a set of sequence of
ordinals of length $\alp$;
\item{(ii)} $\op^\nu \in K$;
\item{(iii)} $\bet < \alp$, $\nu \in Y_\alp$
implies $\nu|\beta \in Y_\bet$,
$\op^{\nu|\bet} \le \op^\nu$ (i.e. $T^\nu
\subseteq T^{\nu|\bet})$, $[i \in T^\nu
\Rightarrow p^{\nu|i} \le p^\nu]$;
\item{(iv)} $Y_0 = \{ <>\}$, $\op^{<>} =
\langle p_i^{<>}: i < \lam\rangle$, \quad
$p_i^{<>} = \emptyset_Q$ (the minimal
element);
\item{(v)} for $\del$ limit, $Y_\del =
\{\eta:\eta$ a sequence of ordinals $\ell
g(\eta) = \del$, $(\forall i < \del)\eta| i
\in Y_i\}$;
\item{(vi)} $\{T^\nu: \nu \in Y_\alp\}$ is
a maximal antichain;
\item{(vii)} for $\alp$ limit, $\nu \in
Y_\alp$, $T^\nu = \bigcap\limits_{\alp' <
\alp} T^{\nu|\alp'}$;
\item{(viii)} for $i \in T^\nu$, $\nu \in
Y_{\alp+1}$, $p_i^\nu \buildrul Q\under
\vline^{''} i \in \su A_{j(\nu)}^{\alp
^{''}}$.

There is no problem to do this for $\alp$  
limit, (vi) is preserved as $\forall \Tet <
\lam[\Tet^\kap < \lam]$ and $\{\del: cf
\del \ge \kap\} \in D$ by an assumption and
Fact 1.3.

Now as $D$ is $\kap$-presaturated there are
$B \in D^+$, $y_\alp \subseteq Y_\alp$,
$|y_\alp| \le \lam$, such that $[\nu \in
Y_\alp - y_\alp \Rightarrow B \cap T^\nu
\not\in D^+]$.  So there is $C \in D$ such
that $(\forall \alp < \kap(*)) (\forall
\nu_1 \ne \nu_2 \in y_\alp)$ $[T^{\nu_1}
\cap T^{\nu_2} \cap C = \emptyset]$.  Let
$B' = \bigcap\limits_{\alp < \kap(*)}
\bigcup\limits_{\nu \in y_\alp} (T^\nu \cap
C)$, then $B' \supseteq B$ mod $D$, hence
$B' \in D^+$.
\hb
Apply (d) and get $\op^*$.

The following proposition shows that it is
possible to remove the assumptions on
cofinality used in Proposition 1.10(1).

\proclaim Proposition 1.11.  Suppose $\kap
< \lam$ are regular cardinals, $Q$ a
$\kap$-complete forcing notion, $\buildrul
Q \under \vline$ ``$\lam$ is regular
cardinal". Let 
$\su D^Q$ be the $Q$-name of the normal filter
on $\lam$ which $D$ generates in $V^Q$.
Assume that the following principle holds:
\hb
$(*)^6_{D,K,Q}$:
\item{(a)} $K$ is a set, its elements are of
the form $\op = \langle p_i: i \in
T\rangle$, $p_i \in Q$, $T = T_\op \in
D^+$;
\item{(b)} if $q \in Q$, $T \in D^+$ then
there is $\op \in K$, $\wedge_i q \le p_i$,
$T \supseteq T_\op$;
\item{(c)} let $\su \tau_\op = \{i: p_i
\in \su G\}$ we assume $p_i \vline \su
\tau_\op \in (D^Q)^+$ (or just $\vline
\hskip -.5truecm /\enskip \su \tau_\op \not\in (\su
D^Q)^+)$;
\item{(d)} if $\op \in K$, $T' \subseteq
T_\op$, $T' \in D^+$ and $ \wedge_{i \in T'} p_i
\le q_i \in Q$ {\it then} for some $T^{''}
\subseteq T'$ and $\otau = \langle r_i : i
\in T^{''}\rangle \in K$ $\wedge_{i \in
T^{''}} q_i \le r_i$;
\item{(e)} for every $q \in Q$, $\su \tau$
s.t. $q \vline \su \tau \in (D^Q)^+$, there
is $\op \in K$, $p_i \ge q$, $q \vline
\tau_\op \subseteq \su \tau$;
\item{(f)} $\langle \emptyset_Q: i < \lam\rangle \in
K$, where $\emptyset_Q$ is the minimal
element of $Q$.

If $D$ is $\kap$-presaturated then $\su
D^Q$ is $\kap$-presaturated.

\subheading{Remark}  We can omit (b) and
get only $\vline \hskip -.5truecm /$\quad ``$\su
D_i^Q | \su T_\op$ is not
$\kap$-presaturated" if $\op \in K$.

\pr Let $\kap(*) < \kap$, $\su I_\alp =
\{\su A_i^\alp: i < \su i_\alp\}$ for $\alp
< \kap(*)$  are
$\kap(*)$ $Q$-names of maximal antichains
of $(\su D^Q)^+$ which form a counterexample
to $\kap$-presaturativity (i.e. at least some
$q_0 \in P$ forces this).

We define by induction on $\alp \le
\kap(*)$  $Y_\alp$ and the function
$j|Y_\alp: Y_\alp \to$ ordinals and $p^\nu
(v \in Y_\alp)$ s.t.
\item{(i)} $Y_\alp$ a set;
\item{(ii)} $\op^v \in K$ let $\op^v =
\langle p_i^v: i \in T^v\rangle, \quad T^v
\in D^+$;
\item{(iii)} for every $v \in Y_\alp$ and
$i \in T^v$ $q_0 \le p_i^v$;
\item{(iv)} $\langle T^v: v \in Y_\alp 
\}$ is a maximal antichain of $D^+$;
\item{(v)} for $v \in Y_\alp$, $i \in T^v$,
$p_i^v \vline^{''}i \in \su
A^{\alp\ ''}_{j(v)}$;
\item{(vi)} if $\eta \in Y_\beta$, $\bet <
\alp$, $v \in Y_\alp$, and $T^v \cap T^\eta
\in D^+$ {\it then} for some $C_{v,\eta}
\in D$, $(\forall i \in T_v \cap T_\eta
\cap C_{v,\eta})$ $[p_i^\eta \le p_i^v]$;
\item{(vii)} $\bet < \alp \Rightarrow
Y_\beta \cap Y_\alp = \emptyset$;
\item{(viii)} for $\beta < \alp$, $v \in
Y_\alp$ $|\{\eta \in Y_\beta: T^\eta \cap
T^v \in D^+\}| \le \lam$.

If $\eta, v \in \bigcup\limits_{\beta <
\alp} Y_\bet$ and
$T^\eta \cap T^v \not\in D^+$, then we
choose $C_{\eta,v} \in D$ disjoint to
$T^\eta \cap T^v$.  (It occurs when $\eta
\ne v \in Y_\bet$.)  

Arriving at $\alp$, let $\{\op^v: v \in
Y_\alp\}$ be maximal such that (i), (ii),
(iii), (v), (vi), (vii), (viii), holds and
$\{T^v: v \in Y_\alp\}$ is an antichain on
$D^+$ (not necessarily maximal).

It suffices to prove that it is a maximal
antichain, so let $T \in D^+$, $T \cap T^v =
\phi \MOD D$ for every $v \in Y_\alp$.  As
$D$ is $\kap$-presaturated, there is $T'
\subseteq T$, $T' \in D^+$ s.t. for every
$\bet < \alp$, $\{\eta \in Y_\bet: T' \cap
T^\eta \in D^+\}$ has cardinality $\le
\lam$, and let it be $\{\eta_{\bet, j}: j <
j_\bet \le \lam\}$.  Let $C = \{\del <
\lam$: if $\bet_1\bet_2 < \alp$, $j_1 <
\del$, $j_2 < \del$, {\it then} $\del \in
C_{\eta_{ \bet_1, j_1}, \eta_{\bet_2,
j_2}}\}$.  Now for each $\bet < \alp$, let
$T'_\beta = \{ \del \in T': (\exists j <
j_\bet \cap \del)[ \del \in T^{\eta \bet,
j}]\}$.  Then $T'_\bet \subseteq T'$ and
$T' - T'_\bet \not\in D^+$ (as $\{T^v: v
\in Y_\bet\}$ is a maximal antichain).  Let
$T^* = C \cap T' \cap \bigcap\limits_{\bet
< \alp} T'_\bet$. So $T' - T^* \not\in
D^+$.  Hence $T^* \in D^+$, $T^* \subseteq
T$.  {\it For every} $\del \in T^*$ and
$\beta < \alp$, there is a unique
$j_\beta^\del < \del$, $\del \in
T^{\eta_{\bet, j^\del_\bet}}$ (if $j_1, j_2$
are candidates use the definition of
$C_{\eta_{ \bet, j_1}, \eta_{\bet, j_2}}$
(there is one as $\del \in T'_\bet$).  Now
for $\bet_1 <\bet_2 < \alp$ (for our
$\del$) by the choice of
$C_{\eta_{\beta_1j_1}, \eta_{\beta_2j_2}}$
$p^{\eta_{\bet_1, j^\del_{\bet_1}}} \le
p^{\eta_{\bet_2, j^\del_{\bet_2}}}$.  Hence
using $\kap$-completeness of $Q$, we can
find $p_\del$, $\wedge_{\bet<\alp} p_\del$
 $\le p_\del$.
Now $\langle p_\del: \del \in T^*\rangle$
satisfies much, almost contradicting the
maximality of $\{\op^v: v \in Y_\alp\}$ and
non-maximality of $\{T^v: v \in Y_\alp\}$
(and the choice of a $j$).  But repairing
this is easy.  Without loss of generality
$T^* \subseteq T^{\eta_{0,0}}$.  Using (d)
of
$(*)^6_{D,K,Q}$ we obtain $\op' = \langle
p'_\del | \del \in T^a\rangle \in K$ s.t.
$T^a \subseteq T^*$ and $\wedge_{i \in
T^a} p_i \le p'_i$. We have to take care of
(v).  So we have $\op' = \langle p'_\del:
\del \in T^a\rangle \in K$ $\wedge_{v \in
Y_\alp} T^a \cap T^v = \emptyset$, $\op'$
satisfies everything except (v).  Now
$\vline\hskip -.5truecm /\enskip\ \tau_{\op'}
\not\in (\su D^+)^Q$ so for some $r\in Q$
$r \vline \tau_{\op'} \in (\su D^+)^Q$.
Hence
for some $r_1$, $r \le r_1 \in Q$, and $j$
the following holds: \hb
$(*)\enskip \tau_1 \vline \tau_{\op'} \cap
\su A_j^\alp \in D^+$.

For each $i \in T^*$ choose if possible
$p^{''}_i$, $p'_i \le p_i^{''} \in Q$
$p^{''}_i \vline^{''} i \in
A_j^{\alp^{''}}$.  Set $T^b = \{i:
p^{''}_i$ is defined$\}$.  It is necessarily
in $D^+$ (otherwise this contradicts
$(*)$).  Applying (d) of $(*)^6_{K,Q,D}$ we
get a final $\op$ contradicting ``$Y_\alp$
maximal but $\{T^v: v \in Y_\alp\}$ is
not".\hb
For $\alp = 0$, use (b) with our
$q_0$.\hfill$\bigsquare$

\proclaim Proposition 1.12.  Let $D,D_1$
are normal filters over a regular cardinal
$\lam$.  Suppose that the following holds
$$C_{D,D_1}: \quad \hbox{for every} \quad T
\in D^+, \quad F(T) = \{\del < \lam: T \cap
\del \ \hbox{is stationary}\} \in D_1^+\
.$$
Let $Q$ be a $\kap$-closed forcing for a
regular $\kap < \lam$.  Then in $V^Q$
$C_{D, D_1}$ is not forced to fail provided
that
$$\eqalign{
(*)^7_{D,K,D_1,Q}\qquad &(\alp) (*)^6_{D,K,Q}
\quad {\rm or} \quad (*)^3_{D,K,Q}\quad {\rm
and}\cr
&(\beta_1)\quad {\rm for} \quad \op \in
K\cr
&\{\del < \lam: \su \tau_\op \cap
\del\quad \hbox{is forced to be
stationary}\} \in D_1^+\cr
&{\rm or}\quad (\bet_2)\quad{\rm for}\quad
\op \in K\cr
&\{\del < \lam: \su \tau_\op \cap \del\quad
\hbox{is not forced to be nonstationary}\}
\in D_1^+\cr
&{\rm and}\quad Q\quad {\rm satisfies}\quad
\lam- {\rm c.c}\ .\cr}$$

\pr Suppose otherwise.  Let $\vline^{''}
\su S \in (\su D^Q)^+$ and $F(\su S)
\not\in (\su D^Q_1)^{+\ ''}$.  Set $T = \{i <
\lam|$ some $p_i$ forces $'' i \in \su
S^{''}\}$.  Using (d) for $\op = \langle
\emptyset: i < \lam\rangle$ (see (f)) find
$\OR \in K$, $\OR = \langle r_i|i \in
T'\rangle$ s.t. $T' \subseteq T$ and $r_i
\ge p_i$.  Then $\vline \su \tau_\OR
\subseteq \su S$.  If $(\bet_1)$ holds,
then we are done.  If $(\beta_2)$ is true,
then by $\lam$-c.c. of $Q$ and normality of
$D_1$ find $C \in D_1$ so that $\vline^{''}
\check C \cap F(\su S) = \emptyset^{''}$.
Then also $\vline^{''} \check C \cap F(\su
\tau_\OR) = \emptyset^{''}$.  The set
$\{\del \in C|\ \hbox{for some}\ q_\del,\ 
q_\del \vline \su \tau_\OR \cap \check
\del$ is stationary$\}$ is in $D_1^+$.
Hence it is nonempty.  So there is $\del
\in C$ s.t. $q_\del \vline \su\tau_\OR \cap
\check \del$ is stationary.  But then
$q_\del \vline^{''} \check \del \in
F(\su\tau_\OR)$ and $\check \del \in \check
C^{''}$ which contradicts the choice of
$C$.\hfill$\bigsquare$

\proclaim Lemma 1.13.  Suppose that $\mu$
is a regular cardinal, $\lam > \mu$ is an
inaccessible $Q = {\rm Col}(\mu, < \lam)$, and
$D$ is a normal filter on $\lam$.  Let $T
\in D^+$ and $\op = \langle p_i|i \in T
\rangle$ be a sequence of conditions in
$Q$.  Then there is $C \in D$ such that for
every $i \in T \cap C$ $p_i \vline \su
\tau_\op \in (\su D^Q)^+$, where $\su
\tau_\op = \{ i: p_i \in \su G\}$.

\pr Set $S = \{i \in T: p_i\vline\hskip
-.5truecm / \enskip\ \tau_\op \in (\su
D^Q)^+ \}$.  If $\lam - S \in D$, then let
$C = \lam - S$.  It will be as required
since
$$\vline\ (\su \tau \subseteq \tau_\op
\quad {\rm and} \quad \tau_\op -  \su
\tau\subseteq \lam - C)$$
i.e $\su \tau$ and $\tau_\op$ are the same
mod $\su D^Q$ where $\tau = \{i \in T \cap
C|p_i \in \su G\}$.

Let us now assume that $S \in D^+$.  For
every $i \in S$ there is $q_i \ge p_i$ $q_i
\vline \su \tau_\op \not\in (\su D^Q)^+$.
Set $\su \tau^* = \{ i \in S: q_i \in \su
G\}$.  Since $Q$ satisfies $\lam$-c.c. and
$D$ is $\lam$-complete, there exists $q$
forcing $''\su \tau^* \in (\su D^Q)^{+\ ''}$.
Then for some $q' \ge q$, $i_0 \in S$ $q'
\ge q_{i_0}$.  So $q' \vline^{''} \{ i \in
S: q_i \in \su G\} \in (D^Q)^{+\ ''}$.  Hence
$q' \vline ^{''} \su \tau_\op \in
(D^Q)^{+\ ''}$.  Which is impossible since
$q_{i_0} \vline^{''} \tau_\op \not\in
(D^Q)^{+\ ''}$.  Contradiction.  So $S
\not\in D^+$.

\subheading{Remark 1.13A} It is possible to
replace the Levy collapse by any
$\lam$-c.c. forcing.

\proclaim Proposition 1.14.  Suppose that
$\mu$ is a regular cardinal $\lam > \mu$ is
an inaccessible, $Q$ is the Levy collapse
Col$(\mu, < \lam)$ and $D, D_1$ are normal
filters over $\kap$.  {\it Then}
\item{1)} $(*)^3_{D,K,Q}$ holds if $\{\del
< \lam| {\rm cf} \del \ge \mu\} \in D$ and
we let $K = \{ \langle p_i|i \in T \rangle
|T \in D^+, \enskip p_i \in Q\}$;
\item{2)} $(*)^6_{D,K,Q}$ holds if we let
$K = \{ \langle p_i | i \in T \cap C
\rangle | T \in D^+$, $p_i \in Q$, $C$ is
as in Lemma 1.13$\}$;
\item{3)} $(*)^4_{D,K,Q}$ holds if $\{\del
< \lam| {\rm cf} \del \ge \mu\} \in D$ and
$K$ is as in 2);
\item{4)} $(*)^7_{D,K,D_1, Q}$ holds if
$\{\del < \lam| {\rm cf} \del < \mu\} \in
D$, $\{\del < \lam|\del$ is an inaccessible
$\} \in D_1$ and $K$ is as in 2).

\pr 1), 2), 3) easy checking.  4).  Suppose
that in $V$ $C_{D,D_1}$ holds.  Let us show
the condition $(\bet_2)$ of
$(*)^7_{D,K,D_1, Q}$.  Let $\op = \langle
p_i | i \in T \rangle \in K$ and $C \in
D_1$.  We should find some $\del \in C$ and
$q \in Q$ $q \vline \su \tau_\op \cap \del$
is stationary.

Using $C_{D,D_1}$ find an inaccessible
$\del \in C$ such that $\del \in C$, $T
\cap \del$ is stationary and for every $i <
\del$ $p_i \in {\rm Col} (\mu, < \del)$.
Now, as in Lemma 1.13, there is a
stationary $S \subseteq T \cap \del$ such
that every $p_i$ $(i \in S)$ forces in Col
$(\mu, < \del)$  
 that
$``\su\tau_\op \cap \del$ is stationary".
If $D$ concentrates on ordinals of
cofinality $<\mu$, then without loss of
generality all elements of $T$ are of some
fixed cofinality $<\mu$.  The forcing
Col$(\mu, < \lam)$ is $\mu$-closed, so it
would preserve the stationarity of
$(\tau_\op \cap \del)^{{\rm Col}
(\mu,<\del)}$.\hfill$\bigsquare$

\proclaim Corollary 1.15.  Let $\lam$ be an
inaccessible, $\mu < \lam$ be a regular
cardinal and $Q = {\rm Col} (\mu,< \lam)$.
Then the following properties are preserved
in the generic extension (i.e. for any
normal filter $D$ on $\lam$)
\item{(1)} $\kap$-presaturatedness, for
$\kap \le \mu$;
\item{(2)} the existence of winning
strategy for II in games defined in
Definition 1.7;
\item{(3)} the nonexistence of winning
strategy for I;
\item{(4)} reflection of stationary subsets
of ordinals of cofinality $<\mu$ provided
that $\lam$ is Mahlo and in $V$ the
reflecting ordinals are inaccessibles.

It is possible to use (4) in order to give
an alternative to the Harrington-Shelah [H-S]
proof of ``every stationary subset of
$\aleph_2$ consisting of ordinals of
cofinality $\ome$ reflects" from a Mahlo
cardinal.

Let $\lam$ be a Mahlo cardinal.  Use the
Backward Easton iteration in order to add
$\alp^+$ Cohen subsets to every regular
$\alp< \lam$.  Now iterate the forcing for
shooting clubs through compliments of
nonreflecting subsets of $\{\alp < \lam|
{\rm cf} \alp = \ome\}$ as it is done in
[H-S].  Such forcing will preserve all the
cardinals.  Since it is possible to pick a
submodel $N$ s.t. $\del = N \cap \lam$ is
an inaccessible cardinal and use Cohen
subsets of $\del$ for the definition of
$N$-generic clubs.  We refer to [S1] for
similar construction with $\diamond$.

Let us now give an example of a forcing
notion satisfying the preservation
conditions but not $\lam$-c.c.  The forcing
notion we are going to consider was
introduced by J. Baumgartner [B] and in a
slightly different form by U. Avraham [A].

\proclaim Proposition 1.16.  Suppose $\mu <
\lam$ are regular cardinals, $\forall \Tet
< \lam$ $(\Tet^{< \mu} < \lam)$, $D$ is the
closed unbounded filter on $\lam$
restricted to some cofinalities $\ge \mu$
(or to a stationary set consisting of
ordinals of such cofinalities) and $S \in
D$ is a stationary subset of $\lam$
containing all ordinals of cofinality
$<\mu$.  Let $Q = \{A|A$ is a set of
pairwise disjoint closed intervals
$\subseteq \lam$ of cardinality $< \mu$,
and $[\alp, \bet] \in A$ implies $\alp \in
S\}$ ordered by inclusion (so $\buildrul Q
\under \vline^{''} S$ contains a club $\su
C_Q = \{\alp |$ for some $\bet$ $[\alp,
\bet] \in \su G_Q^{''})$.  Define $K = \{
\langle p_i | i \in T \rangle | T \in
D^+,\enskip T \subseteq \{\del \in S| {\rm
cf} \del \ge \mu \}$, $p_i \in Q$ for some
$j \ge i$, \enskip $[i,j] \in p_i$ and $p_i
\vline \tau_\op \in (D^Q)^+$ for every $i
\in T\}$.
\hb
Then $(*)^6_{D,K,Q}$, $(*)^3_{D,K,Q}$,
$(*)^4_{D,K,Q}$ hold.

\pr
Let $\op = \langle p_i | i \in T \rangle$
be so that $T \in D^+$, $T \subseteq \{
\del \in S | {\rm cf} \del \ge \mu\}$, $p_i
\in Q$ and for some $j \ge i$, $[i,j] \in
p_i$.  Let us show that $\vline\hskip
-.5truecm / \enskip$ ``$\su \tau_\op$
is not stationary".  Suppose otherwise.
Let $\su E$ be a name of a club such that
$\vline \su E \cap \su \tau_\op =
\emptyset$.  For every $i \in T$ pick $q_i
\ge p_i$ forcing $''i \not\in \su E^{''}$.
Pick an increasing continuous sequence
$\langle M_i | i < \lam\rangle$ of
elementary submodels of some $H(\tau)$ for
$\tau$ big enough so that $|M_i| < \lam$
and $M_i \cap \lam \in \lam$.  Set $C^* =
\{ \del| M_\del \cap \lam = \del\}$.
Clearly, $C^*$ is a club.  Pick $\del^* \in
C^* \cap T$.  Then $q_{\del^*} \vline\del^*
\in \su E$ since otherwise for some $\alp <
\del^*$, $r \ge q_{\del^*}$ $r \vline \su E
\cap \del^* \subseteq \alp$.  But since cf
$\del^* \ge \mu$ and $|r| < \mu$, for some
$\bet, \alp < \bet < \del^*$ there exists
$r' \in Q \cap M_\bet$, $r' \ge r| \del^*$
forcing ``$\su E \cap \del^* \not\subseteq
\alp$".  This leads to the contradiction
since $r' \cup r^* \in Q$.

So for every $p$ as above some $q \in Q$
forces ``$\su \tau_\op$ is stationary".
Now the arguments of Lemma 1.13 apply.  So
for every $\op = \langle p_i | i \in T
\rangle$ as above there exists $C \in D$
such that for every $i \in T \cap C$ $p_i
\vline \su \tau_\op \in (\su D^Q)^+$.

Let us show that for any $q \in Q$, $\su T$
s.t. $q \vline \su T \in (\su D^Q)^+$ there
exists $\op \in K$, $p_i \ge q$ and $q
\vline ^{''} \tau_{\op} \subseteq \su T
\MOD (\su D^Q)^{+\ ''}$.
The set $T = \{ \del|\hbox{for some}\
p_\del \ge q,\ p_\del \vline \del \in \su T
\cap \su C_Q\}$ is in $D^+$.  But if
$p_\del \vline \del \in \su C_Q$, then for
some $\del'[\del, \del'] \in p_\del$.  Now
fix $\langle p_\del| \del \in T \rangle$ as
above in order to obtain $\op \in K$ as
required.

The checking of the rest of the conditions
is routine.\hfill$\bigsquare$

What happens if $D$ concentrates on small
cofinality?  Does forcing with $Q$ of
Proposition 1.6 preserve $<
\mu$-presaturedness?  Strengthening the
assumptions it is possible to obtain a
positive answer.  Namely the following
holds.

\proclaim Proposition 1.17.  Let $\mu,
\lam, S, Q$ be as in Proposition 1.16,
assume that $D$ is a club filter restricted
to some cofinality $< \mu$ (or to a
stationary set of such cofinality) and
there exists a set $S^- \in D$ consisting of
ordinals of cofinality $< \mu$ and for
every $\del \in S^-$ there is a set $A_\del
\subseteq \del$, $|A_\del| < \mu$
consisting of ordinals of cofinality $\ge
\mu$, and the following holds:\hb
$(*)$ ``for every club $C \subseteq \lam$
$\{ \del \in S^-|\sup(A_\del \cap C) =
\del\} \in D$".\hb
Let $K = \{ \langle p_i|i \in T \rangle | T
\in D^+,\enskip p_i \in Q,T \subseteq S^-$,
for every $i \in T (\exists \alp < i)
(\forall \xi \in A_i - \alp) (\exists \eta)
([\xi, \eta] \in p_i)\}$. Then
$(*)^6_{D|\su S^*, K, Q}$ holds, where $\su
S^* = \{ \del \in S^-|\sup (A_\del - \su
C_Q) < \del\}$.

\subheading{Remarks}
\item{(1)} The property $(*)$ is true in
$L$ and also it can be easily forced.
\item{(2)} If we are interested only in a
condition forcing $< \mu$-presaturation
then there is no need in $\su S^*$.

\pr Let us check that for $\op = \langle
p_i|i \in T\rangle \in K\enskip \vline\hskip
-.5truecm / \enskip``\tau_\op$ is not stationary".
Suppose otherwise, let $\su C$ be a name of
club disjoint to $\su \tau_\op$.  Pick $q_i
\ge p_i$ forcing ``$i \not \in \su C$".
Let $\langle M_\alp| \alp < \lam\rangle$ and
$C^*$ be as in the previous proposition.
Pick $\del^* \in C^* \cap T$ such that
otp$(C^* \cap \del^*) = \del^*$ and sup
$(A_{\del^*} \cap C^*) = \del^*$.  By the
assumption $q_{\del^*} \vline^{''} \del^*
\not\in \su C^{''}$.  So for some $r \ge
q_{\del^*}$ and $\alp < \del^*$ $r
\vline^{''} \su C \cap \del^* \subseteq
\alp^{''}$.  But it is possible to find $\del
\in (A_{\del^*} \cap C^*) - \alp$ such that
$r| \del \in M_\del$.  So some $r' \ge r|
\del$, $r' \in M_\del$ forces ``$\su C \cap
\del^* \not\subseteq \alp$".  But $r \cup
r' \in Q$.
Contradiction.\hfill$\bigsquare$

\subheading{2.~~Constructions of 
Cardinal Preserving Ideals}

By Jech-Magidor-Mitchell-Prikry [J-M-Mi-P]
it is possible to construct a model with a
normal $\mu$-preserving ideal over $\mu^+$
for any regular $\mu$ from one measurable.
Actually $NS^\mu_{\mu^+}$ can be such
ideal.  We shall examine here
$NS_{\mu^+}^r$ for $\tau < \mu$ and
$NS_{\mu^+}$.

In order to formulate the results we need
the following definition of [G4]:

\subheading{Definition 2.0}  Let $\aF = <
\calF(\alp, \bet) | \bet < \gam)$ be a sequence
of ultrafilters over $\alp$.  Let $\del <
\gam$, $\rho > 0$ be ordinals and $\lam \le \alp$ be
a regular cardinal.  Then
\item{(a)} $\del$ is an up-repeat point for
$\aF$ if for every $A \in \calF(\alp, \del)$
there is $\del'$, $\del < \del'<\gam$ such
that $A \in \calF(\alp, \del')$;
\item{(b)} $\del$ is a $(\lam, \rho)$-repeat
point for $\aF$ if $(a1)$ $cf\del = \lam$
and $(a2)$ for every $A \in \cap\{\calF(\alp,
\bet)| \del \le \bet < \del + \rho\}$ there
are unboundedly many $\xi$'s in $\del$ such
that 
$A \in \cap\{\calF(\alp, \xi')| \xi \le \xi' <
\xi + \rho\}$.

If $\aF$ is the maximal sequence of measures
of the core model we will simply omit it.
\proclaim Theorem 2.1.  The exact strength
of
 \item{(1)} $``NS_{\mu^+}$ is
$\mu$-preserving for a regular $\mu >
\aleph_2 + GCH$" or;
\item{(2)} ``$NS_{\mu^+}^{\aleph_0}$ is
$\mu$-preserving for a regular $\mu >
\aleph_2 + GCH$" or;
\item{(3)} ``$NS_{\mu^+}$ is
$\aleph_2$-preserving for a regular $\mu >
\aleph_2 +GCH$" or;
\item{(4)} ``$NS^r_{\mu^+}$ is
$\aleph_2$-preserving for a regular $\mu >
\aleph_2$, $\tau < \mu + GCH$"\hb
is the existence of an up-repeat point.

\pr Use the model of [G4] with a
presaturated $NS_\kap$ over an inaccessible
$\kap$ and the Levy collapse
Col$(\mu,<\kap)$.  By the results of
Section1, $NS_\kap$ will be
$\mu$-preserving in the generic
extension.\hfill$\bigsquare$

The rest of the section will be devoted to
the construction of $NS_\kap$
$\ome_1$-preserving from an $(\ome,
\kap^+ + 1)$-repeat point. The following
theorem will then follow.
\proclaim Theorem 2.2.
\item{(1)} the strength of
``$NS_\kap^{\aleph_0}$ is
$\ome_1$-preserving $+\kap$ is an
inaccessible" is an $(\ome, \kap^+ +
1)$-repeat point;
\item{(2)} the existence of an $(\ome, \kap^+
+ 1)$-repeat point is sufficient for
``$NS_\kap$ is $\ome_1$-preserving $+\kap$
is an inaccessible $+GCH$";
\item{(3)} the strength of
``$NS_{\mu^+}^{\aleph_0}$ is
$\omega_1$-preserving $+GCH$" for a regular
$\mu > \aleph_2$ is $(\omega, \mu)$-repeat
point;
\item{(4)} the existence of an $(\ome,\mu +
1)$-repeat point is sufficient for
$``NS_{\mu^+}$ is $\ome_1$-preserving
$+GCH"$, where $\mu$ is a regular cardinal.

Suppose that $\vec U$ is a
coherent sequence of ultrafilters with
an $(\ome, \kap^+ + 1)$-repeat point $\alp$ at
$\kap$.  Define the iteration $\calP_\tau$ for
$\tau$ in the closure of $\{\beta|(\bet =
\kap)$ or $(\bet < \kap$ and $\bet$ is an
inaccessible or $\bet = \gam + 1$ and
$\gam$ is an inaccessible)$\}$.

On the limit stages use the limit of [G1].
For the benefit of
the reader let us give a precise definition.

Let $A$  be a set consisting of $\alp$'s such
that $\alp <\kap$  and $\alp >0(\alp)>0$.
Denote by $A^\ell$  the closure of the set $\{\alp
+1\mid \alp\in A\}\cup A$.  For every $\alp\in
A^\ell$  define by induction $\calP_\alp$ to be
the set of all elements $p$  of the form
$\langle\su p_\alp\mid\gam\in g
\rangle$  where 
\item{(1)} $g$ is a subset of $\alp\cap A$.
\item{(2)} $g$ has an Easton support, i.e. for
every inaccessible $\bet\le\alp$,  $\bet
>|\dom g\cap\bet |$; 
\item{(3)} for every $\gam\in\dom g$
$p\rhookup\gam =\langle \su p_\bet
\mid \bet\in\gam\cap g\rangle\in\calP_\gam$  and
$p\rhookup\gam\buildrul{\calP_\gam}\under
\vline$
``$\su p_\gam\in
\su Q_\gam$".

Let $p=\langle \su p_\gam
\mid\gam\in g\rangle\ ,\ q=\langle 
\su q_\gam\mid\gam\in f\rangle$ be elements
of $\calP_\alp$.  Then $p\ge q$  ($p$  is
stronger than $q$) if the following holds: 
\item{(1)} $g\supseteq f$
\item{(2)} for every $\gam\in f$ 
$p\rhookup \gam\buildrul{\calP_\gam}\under\vline$
``$\su p_\gam\ge
\su q_\gam$  in the forcing
$\su Q_\gam$"
\item{(3)} there exists a finite subset $b$ of
$f$ so that for every $\gam\in f\bks b$,
$p\rhookup \gam\vline$ ``$
\su p_\gam\ge^* \su q_\gam$
in $\su Q_\gam$".

If $b=\emptyset$,  then let $p\ge^*q$.

Suppose that $\tau$ is an inaccessible and
$\calP_\tau$ is defined.  Define
$\calP_{\tau + 1}$.  Let $C(\tau^+)$ be the
forcing for adding $\tau^+$-Cohen subsets
to $\tau$, i.e. $\{f \in V^{\calP_\tau}|f$
is a partial function from $\alp^+ \times
\alp$ into $\alp$, $|f|^{V^{\calP_\alp}} <
\alp$ and for every $\beta < \alp^+ \{
\beta'| (\bet, \bet') \in {\rm dom}\ f\}$
is an ordinal$\}$.  $\calP_{\tau + 1}$ will
be $\calP_\tau * C(\tau^+) * \calP(\tau,
O^{\vec U} (\tau))$, where
$\calP(\tau,O^{\aU} (\tau))$ is the forcing
of [G1,2] with the slight change described
below.

The change is in the definition of
$U(\tau,\gam,t)$ the ultrafilter extending
$U(\tau,\gam)$ for $\gam < O^{\aU}(\tau)$ and
a coherent sequence $t$ or more precisely
in the definition of the master conditions
sequence.  Let $j_\bet^\tau: V \to
N_\bet^\tau \cong V^\tau/ U(\tau, \beta)$
for $\bet < O^{\aU}(\tau)$.  Pick some well
ordering $W$ of $V_\lam$, for a big enough
$\lam$ so that for every inaccessible
$\del < \lam$, $W|V_\del: V_\del
\leftrightarrow \del$.  Let $\gam$ be some
fixed ordinal below $O^{\aU} (\tau)$.  Let us
for a while drop the indices $\tau, \gam$
in $j_\gam^\tau$, $N_\gam^\tau$.

Let $< \su D_{\gam'} | \gam' < \tau^+ >$ be
the $j(W)$-least enumeration of all
$E$-dense open subsets of $j(\calP_\tau *
C(\tau^+))$ which are in $N$.  Where a
subset $D$ of a forcing notion $\calP$ 
$E$-dense is if for every $p \in \calP$ there
exists $q \in D$ which is an Easton
extension of $p$.  Define an $E$-increasing
sequence $<\su p_{\gam'}|\gam' < \tau^+>$ of
elements of $\calP_{j(\tau)} *
C(j(\tau)^+)/\calP_{\tau + 1}$ so that for
every $\gam' < \tau^+ \su p_{\gam'}$ there
will
be a $j(W)$-least $E$-extension of
$<\su p_{\gam^{''}}| \gam^{''} < \gam'>$ in
$\su D_{\gam'}$ compatible with $j^{''} (G
\cap C(\tau^+))$.

Now as in [G1,2] set $A \in U(\tau, \gam,
t)$ if for some $r$ in the generic
subset of $\calP_\tau * C(\tau^+)$, some
$\gam'<\tau^+$, a name $\su A$ of $A$ and a
$\calP_\tau * C(\tau)$ -- name $\su T$, in
$N$
$$r \cap \{ < \check t, \quad \su T > \}
\cup \su p_{\gam'} \vline \check \tau \in j(\su
A)\ .$$
Note that the set $D = \{ q|q||(\check \tau
\in j(\su A)\}$ is $E$-dense and it belongs
to $N$.  So it appears in the list $<\su
D_{\gam'} | \gam' < \tau^+>$.  Hence some
$\su p_{\gam'}$ is in $D$.

Let $G$ be a generic subset of $\calP_\kap
* C(\kap^+)$.  Recall that by [G2] the
sequence $< U(\gam, \gam', \emptyset) |
\gam' < O^{\aU} (\gam) > $ is commutative
Rudin-Keisler increasing sequence of
ultrafilters, for every $\gam$.

Now let $j^*: V[G]\to M^* \cong V[G]^\kap/U
(\kap, \gam, \emptyset)$.  Then $M^* = M[j^*
(G)]$ for a model $M$ of $ZFC$ which is
contained in $M^*$.  We would like to have
the exact description $M$.

The next definition is based on the Mitchel
notion of complete iteration see [Mi
1,2,3].

\subheading{Definition 2.3}  Suppose that
$N$ is a model of set theory, $\aV$ a
coherent sequence in $N$, $\kap$ is a
cardinal.  The complete iteration $j: N \to
M$ of $\aV$ at $\kap$ will be the direct
limit of $j_\nu: N \to M_\nu$, $\nu <
\ell_\kap$ for some ordinal $\ell_\kap$,
where $\ell_\kap$, $j_\nu$, $M_\nu$ are
defined as follows.  Set $M_0 = N$, $j_0 =
id$, $\aV_0 = \aV$, $C_0(\alp, \bet) =
\emptyset$ for all $\alp$ and $\bet$.  If
$O^{\aV} (\kap) = 0$, then $\ell_\kap = 1$.
Suppose otherwise:
\subheading{Case 1}  $O^{\aV} (\kap)$ is a
limit ordinal.

If $j_\nu$, $M_\nu$, $\aV_\nu$, $\aC_\nu$
are defined then set
$$j_{\nu\nu + 1}: M_\nu \to M_{\nu + 1}
\cong M^{\alp_\nu}/V_\nu (\alp_\nu,
\beta_\nu)\ ,$$
where $\alp_\nu$ is the minimal ordinal
$\alp$ so that
\item{(i)} $\alp \ge \kap$;
\item{(ii)} $\alp$ is less than the first
$\kap'>\kap$ with $O^{\aV} (\kap') > 0$;
\item{(iii)} for some $\bet < O^{\aV_\nu}
(\alp)$ $\aC_\nu (\alp, \bet)$ is bounded
in $\alp$;\hb
and $\bet_\nu$ is the minimal $\bet$
satisfying (iii) for $\alp_\nu$.

If there is not such an $\alp_\nu$ then set
$\ell_\kap = \nu$, $j = j_\nu$, $M =
M_\nu$, $\aC = \aC_\nu$.  Define
$\aC_{\nu+1} | (\alp_\nu, \bet_\nu) =
\aC_\nu | (\alp_\nu, \bet_\nu)$,
$\alp_{\nu+1} = j_{\nu\nu + 1} (\alp_\nu)$
$$\aC_{\nu + 1} (\alp_{\nu + 1}, j_{\nu\nu
+ 1} (\bet)) = \cases{
\aC_\nu(\alp_\nu, \bet_\nu), &if $\bet \ne
\bet_\nu$\cr
\aC_\nu (\alp_\nu, \bet_\nu) \cup
\{\alp_\nu\}, &if $\bet = \bet_\nu$\cr}$$
\subheading{Case 2}  $O^{\aV} (\kap) = \tau
+ 1$.

Define $j_\nu$, $M_\nu$, $\aV_\nu$ and
$\aC_\nu$ as in Case 1, only in (ii) let
$\alp$ be less than the image of $\kap$
under the embedding of $N$ by $V(\kap,
\tau)$.

\proclaim Theorem 2.4.  Let $j^*: V[G] \to
M^*$ be
\item{(a)} the ultrapower of $V[G]$ by the
ultrafilter $U(\tau, \gam, \emptyset)$ \hb
or
\item{(b)} the direct limit of the
ultrapowers with the Rudin-Kiesler
increasing sequence of the ultrafilters
$<U(\tau, \gam', \emptyset)| \gam' < \gam>$.
\hb
Then $M^*$ is a generic extension of $M$,
where $M$ is the complete iteration of
$\aU|(\tau, \gam + 1)$ at $\tau$, if (a)
holds, and of $\aU|(\tau, \gam)$ at $\tau$,
if (b) holds.

\proclaim Main Lemma.  Let $M$ be as in the
theorem and let $i$ be the canonical
embedding of $V$ into $M$.  Then there is
$G^* \subseteq i (\calP_\kap * C(\kap^+))$
so that $G^* \in V[G]$ and $G^*$ is $M[G]$
generic.

\pr Let us prove the lemma by induction on
the pairs $(\tau, \gam)$ ordered
lexicographically.  Suppose that it holds
for all $(\tau, \gam') < (\kap, \alp)$.
Let us prove the lemma for $(\kap, \gam)$.

\subheading{Case 1}  $\alp = \alp' + 1$.

Let $N$ be the ultrapower of $V$ by
$U(\kap, \alp)$ and $j:V \to N$ the
canonical embedding.  We just simplify the
previous notations, where $N =
N^\alp_\kap$, $j = j^\alp_\kap$.  Then $M$
is the complete iteration of $j(\aU)| \kap
+ 1$ at $\kap$ in $N$ if $\alp'$ is a
limit ordinal, or for successor $\alp'$, the
above iteration should be performed
$\ome$-times in order to obtain $M$.  Let
us concentrate on the first case; similar
and slightly simpler arguments work for
the second one.

Denote by $k: N \to M$ the above iteration.
Then the following diagram is commutative
$$\matrix{
&&N&\cr
&j&&\cr
&\nearrow&&\cr
V&&\Bigg\downarrow&k\cr
&i&&\cr
&\searrow&&\cr
&&M&\cr
\cr}$$

By the inductive assumption there exists
$G' \in N[G]$ $M$-generic subset of
$\calP_{k(\kap)} * C((k(\kap))^+)$. 

If $\alp' = 0$ then $k(\alp') = 0^{i(\aU)}
(k(\kap)) = 0$ and $C(k(\kap)^+)$ is only
the forcing used over $k(\kap)$.
Suppose now that $\alp'>0$.  Then $k(\alp')
> 0$ and the forcing over $k(\kap)$ is
$C(k(\kap)^+)$ followed by $\calP(k(\kap),
k(\alp'))$.  Let us define in $N[G]$ a
$M[G']$-generic subset of $\calP(k(\kap),
k(\alp'))$. Recall that a generic subset of
$\calP(k(\kap), k'(\alp))$ can be
reconstructed from a generic sequence
$b_{k(\kap)}$ to $k(\kap)$, where 
$b_{k(\kap)}$ is a combination of a cofinal
in $k(\kap)$ sequences so that
$b^{-1}_{k(\kap)} (\{n\})$ is a sequence
appropriate for the ultrafilters $i(U)
(k(\kap), \del)$ with $\del$ on the depth
$n$.  We refer to [G2] for detailed
definitions.

Let us use the indiscernibles of the
complete iteration $\aC$ in order to define
$b_{k(\kap)}$.  Namely, only $<\aC(k(\kap),
\del)| \del < k(\alp') >$ will be used.

Set $b_{k(\kap)} = \{ < \tau, n, \xi > |
n<\ome, \xi < k(\kap)^+$, $\tau<k(\kap)$,
for some $\del < k(\alp')$ coded by
$<n,\xi>$ in sense of [G2] $\tau \in \aC
(k(\kap), \del))\}$.

Let $G^{''}$ be the subset of
$\calP(k(\kap), k(\alp'))$ generated by
$b_{k(\kap)}$.  Let us show that $G^{''}$ is
$M[G']$-generic subset.  Let $D \in M[G']$
be a dense open subset of $\calP(k(\kap),
k(\alp'))$ and $\su D$ its canonical name.
Since $M$ is the direct limit for some
$\nu$ less than the length of the complete
iteration, for some $\su D_\nu \in M_\nu$,
\enskip $\su D$ is the image of $\su
D_\nu$.  Let us work in $M_\nu$.  Denote by
$G_\nu$ the appropriate part of $G'$ and by
$t$ a coherent sequence which generates
$b_{k(\kap)} | \alp_\nu$.  Let $k_\nu: M \to
M_\nu$ be the part of the complete
iteration $k$ on the step $\nu$.

\subheading{ Case 1.1} $cf^{M_\nu} (k_\nu
(\alp')) < k_\nu (\alp')$.

Then $k_\nu(\alp')$ changes its cofinality
to $cf^{M_\nu} (k_\nu (\alp'))$ after the
forcing with $\calP(k_\nu (\kap), k_\nu
(\alp'))$.  

For every $T$ so that $<t, T> \in
\calP(k_\nu (\kap), k_\nu (\alp'))$, there
exists $i < cf^{M_\nu} (k_\nu (\alp'))$ and
$T^*$ so that $< t, T^*>$ is a condition
stronger than $<t, T>$ and $<t,T^*>$ forces
``some $<t',T'> \in \check D_\nu$ with $t'$
on the level $i$ is in the generic set".
In order to find such $T^*$ just use the
$k_\nu(\alp')$-completeness of the
ultrafilters involved in the forcing
$\calP(k_\nu(\kap), k_\nu (\alp'))$ and the
Prikry property.  Now for every $t' \in
T^*$ which is appropriate for $i$ or some
$j \ge i$ there exists $T'$ such that $<t',
T'> \in D_\nu$.  The same property remains
true for $k_{\nu'} (T^*)$ for every $\nu
\le \nu' <$ length of the iteration $k$.
Pick $\nu' \ge \nu$ to be large enough in
order to contain elements of $b_{k(\kap)}$
appropriate for $i$.  Let $t'$ be a
coherent sequence generating $b_{k(\kap)} |
\alp_{\nu'}$.  It is possible to pick such
$t'$ in $k_{\nu'} (T^*)$.  But then for
some $T'$ $<t', T' \cap k_{\nu'} (T^*) >
\in k_{\nu \nu'} (D_\nu)$.  The image of $<
t', T' \cap k_{\nu'} (T^*) >$ under the
rest of the iteration will be in $G^{''}$.

\subheading{Case 1.2}   
$cf^{M_\nu} (k_\nu (\alp')) = k_\nu
(\alp')$.

Then $k_\nu (\alp')$ changes its cofinality
to $\ome$.

For every $T$ so that $< t, T> \in
\calP(k_\nu (\kap), k_\nu(\alp'))$ there
exist $n < \ome$ and $T^*$ so that $< t,
T^*>$ is a condition stronger than $<t,T>$
and $< t, T^*>$ forces ``some $<t',T'> \in
\check D_\nu$, with $t'$ containing the
first $\check n$ elements of the canonical
$\ome$-sequence to $k_\nu (\kap)$, is in
the generic set".

Pick  $\nu' \ge \nu$ in order to first reach
$n$ elements of the canonical
$\ome$-sequence to $k(\kap)$.  Then proceed
as in Case 1.1.

\subheading{Case 1.3}  $cf^{M_\nu} (k_\nu
(\alp')) = (k_\nu (\alp'))^+$ in $M_\nu$.

Let $D$ be a dense open subset of
$\calP(k_\nu(\kap), k_\nu(\alp'))$.  Set
$D' = \{ < p,  T> \in \calP(k_\nu(\kap),
k_\nu(\alp'))|$ for some level $\del <
(k_\nu (\kap))^+$, for every $< \xi_1 \nek
\xi_n>$ s.t. $\xi_n \in Suc_{T,\mu}(p^\cap
< \xi_1 \nek \xi_{n - 1})$ for some $\mu \ge
\del$, $< p^\cap < \xi_1 \nek \xi_n >$,
$T_{p^\cap < \xi_1\nek \xi_n>} > \in D\}$.

\proclaim Claim. $D'$ is a dense open.

The proof is similar to Lemma 3.11 of [G1].
Define $D_1 = D'$, for every $n < \ome$ set
$D_{n+1} = D'_n$ and finally $D_\ome =
\cap_{n < \ome} D_n$.

\proclaim Claim. For every condition $<
p, T>$ there exists a stronger condition $<
p, T^*>$ in $D_\ome$.

We refer to Lemma 1.4 [G1] for the proof.
Just replace $\sig$ there by ``$\in
D_\ome$".

Let $<t, T^*> \in D_\ome$.  Then for some $n
$ $<t, T^*> \in D_n$.  Now by simple
induction it is possible to show that there
is $\del < (k_\nu(\kap))^+$ so that for
every $\xi \in Suc_{T^*,\del} (t) < t^\cap
< \xi >$, $T^*_{t^\cap < \xi> } > \in
D$.

Now pick a large enough part of the
iteration $k$ to reach the level $\del$.
Continue as in Case 1.1.

It completes the definition of a generic
subset of $\calP_{k(\kap) + 1}$.  Let us
refer to it as $G_{k(\kap) + 1}$.  We now turn
to the construction of the generic
object for the forcing between $k(\kap) +
1$ and $i(\kap) + 1$.  Let us define it by
induction on $\del$, $k(\kap) + 1 \le \del
\le i(\kap) + 1$.  Suppose that for every
$\del' < \del$ a $M$-generic subset
$G_{\del'}$ of $\calP_{\del'}$ is defined
in $V[G]$.  Define $G_\del$.

If there is some $\tau > \del$ and an
indiscernible $\tau'$ for it 
$\tau' \le \del$, then use the inductive
assumptions to produce $G_\del$.  Suppose
now that there is no $\tau, \tau'$ as
above.  Notice, that then $\del =
k(\del^*)$ for some $\del^* \le \del$.  Let
us split the proof according to the
following two cases.
\subheading{Case A} $\del = \del' + 1$.

Then $\calP_\del/G_{\del'} = C(\del^+) *
\calP(\del, 0^{k(\aU)} (\del))$ if
$0^{k(\aU)} (\del) > 0$ and
$\calP/G_{\del'} = C(\del^+)$ otherwise.
So $G_\del$ will be $G_{\del'} * (G' *
G^{''})$, where $G^{''}$ may be empty.

Let us define first $G' \subseteq C(\del^+)$.
We use as inductive assumption that
$k(p_\nu)\mid \del \in G_{\del'}$ for every
$\nu < \kap^+$ where $< p_\nu|\nu <
\kap^+>$ is the master condition sequence
for $U(\kap, \alp, \emptyset)$.  Set $G' =
\cup \{k(p_\nu) (\del) | C(\del^+) | \nu <
\kap^+\}$.

Let us check that $G'$ is
$M[G_{\del'}]$-generic subset of
$C(\del^+)$.  Suppose that $D \in
M[G_{\del'}]$ is a dense open subset of
$C(\del^+)$.  Let $\su D$ be a canonical
$P_{\del'}$-name of $D$.  By the assumption
on $\del$ there are indiscernibles $\kap <
\alp_{\nu_1} < \cdots < \alp_{\nu_n}< \del$
such that the support of $\su D$ is a
subset of $\{ \kap, \alp_{\nu_1} \nek
\alp_{\nu_n}\}$.  Since the forcing
$C(\del^+)$ is $\del$-closed $\su D$ can be
replaced by it dense open subset with
support $\kap$ alone.  Let us assume that
$\su D$ is already such a subset.  Then
$\su D = k(\su D^*)$ for dense open $\su
D^* \in N$ subset of $C((\del^*)^+)$.  By
the choice of the master condition sequence
for some $\nu < \kap^+$
$$p_\nu | \del^* \vline p_\nu (\del^*) |
C((\del^*)^+) \in \su D^*\ .$$
But this implies
$$k(p_\nu) (\del) | C(\del^+) \in D\ .$$
So $G' \cap D \ne \emptyset$.

Suppose now that $0^{k(\aU)} (\del) > 0$.
We need to define $G^{''}$ a $M[G_{\del'} *
G']$-generic subset of $\calP(\del,
0^{k(\aU)} (\del))$.  In this case 
$\del$ is a limit of indiscernibles, i.e.
$\aC(\del, \tau)$ is unbounded in $\del$
for every $\tau < 0^{k(\aU)} (\del ))$.  So
we are in the situation considered above.
The only difference is that some $p'_\nu s$ 
may contain information about the generic
sequence $b_\del$ to $\del$.  In order to
preserve $k(p_\nu)| \del + 1$ in the
generic set, we need to start $b_\del$
according to $k(p_\nu) (\del)$.  Notice,
that further elements of the master
condition sequence do not increase the
coherent sequence given by  $p_\nu$.

\subheading{Case B}  $\del$ is a limit
ordinal.

Set $G_\del = \{ p \in \calP_\del |$ for
every $\del' < \del\enskip p| \del' \in
G_{\del'}\}$.  Let us show that $G_\del$ is
$M$-generic subset of $\calP_\del$.
Consider two cases.

\subheading{Case B.1}  There are
unboundedly many in $\del$ indiscernibles
for ordinals $\ge \del$.

Since $\del = k(\del^*)$, $\del$ is a limit
of indiscernibles for $\del$.  Then $\del$
is measurable in $M$ and the direct limit
is used on the stage $\del$.  Let $< \tau_i
| i < \lam>$ be a cofinal sequence of
indiscernibles.

Let $D \in M$ be a dense open subset of
$\calP_\del$.  Then $\{\tau < \del | D \cap
H(\tau)$ is a dense open subset of
$\bigcup\limits_{\tau' < \tau}
\calP_{\tau'}\}$ contains a club in $M$.
The sequence $< \tau_i | i < \lam>$ is
almost contained in every club of $M$ so
there is $i_0$ s.t. $D \cap \calP_{\tau_0}$
is dense.  But then $G_{\tau_0} \cap D \ne
\emptyset$.  Hence $G_\del \cap D \ne
\emptyset$.

\subheading{Case B.2} There are only
boundedly many in $\del$ indiscernibles
for ordinals $\ge \del$.

Then, since $\del = k(\del^*)$, there is no
indiscernibles for ordinals $\ge \del$
below $\del$.  So there are unboundedly
many in $\del$ ordinals $\tau$ which are in
the range of $k$.

Let $D$ be a dense open subset of
$\calP_\del$ in $M$.  Define $D' = \{ p \in
\calP_\del|$ for some $\bet < \del\enskip p|
\beta \vline ^{''} p\backslash \beta \in D /
\su G_\bet ^{''}\}$.  Set $D_0 = D$, $D_{n+1} =
D'_n$ for every $n < \ome$ and let $D_\ome
= \bigcup\limits_{n < \ome} D_n$.

\proclaim Claim B.2.1.  $D_\ome$ is
$E$-dense subset of $\calP_\del$.

\pr Let $p \in \calP_\del$ define $p^*_E
\ge p$ as in Lemma 1.4 [G1] where $\sig$ is
replaced by belonging to $D_\ome$.  Suppose
that $p^* \not\in D_\ome$.  Then there is
$p' \in D_\ome$, $p' \ge p^*$.  Let $\bet
\in $ domp$^* \cap$ domp$'$ be the last on
which an information about cofinal
sequences is added.  But then for some
$p^{''} \in \calP_\bet$ $p^{''\cap} p^*
(\bet) \vline p^* \backslash \bet \in D_1 /
G_\bet$.  Then, as in Lemma 1.4 [G1], it is
possible to go down until finally for some
$n$ and some $p^{''} \in
\calP_{\min (domp^*)}$
$$p^{''} \vline p^* \in D_n/ G_{\min\
(domp^*)}$$
which contradicts the definition of $p^*$.

$\bigsquare$ of the claim.

It is enough to show that $G_\del \cap
D_\ome \ne \emptyset$.  Since then, for
some $n$, $G_\del \cap D_n \ne \emptyset$.
Now use the fact that all initial segments of
$G_\del$ are $M$-generic.  So let us prove
for every $E$-dense set $D$ that $G_\del
\cap D \ne \emptyset$. 
Without loss of generality we can assume
that $D$ is in the range of $k$,  since it
is
always  possible to find some $\tau <
\del$ in the range of $k$ above the support
of $D$ and the intersection of $\tau$
$E$-dense open subsets of
$\calP_\del/G_\tau$ is $E$-dense open.

Let $D^* \in N$ be an $E$-dense open subset
of $\calP_{\del^*}$ so that $k(D^*) = D$.
By the definition of the master condition
sequence $< p_\nu | \nu < \kap^+>$, for
some $\tau_0 < \kap^+$ $p_{\nu_0} | \del^*
\in D^*$.  Then $k(p_{\nu_0} )| \del \in
D$.  But, by the choice of $<G_{\del'}
|\del' < \del>$ $k(p_{\nu_0})| \del' \in
G_{\del'}$ for every $\del' < \del$.  So
$k(p_{\nu_0})| \del \in G_\del$.

So $G_\del$ is an  $M$-generic subset of
$\calP_\del$.  It completes Case B and
hence also Case 1 of the lemma.

\subheading{Case 2}  $\alp$ is a limit
ordinal.

Use the inductive assumption and the
definitions of the generic sets of Case 1.
Define a generic subset of
$\calP_{k(\kap)}$ as in Case B.1.

$\bigsquare$ of the lemma.

Let $G^* \in V[G]$ be the $M$-generic
subset of $\calP_{i(\kap)} * C(i(\kap)^+)$
defined in the Main Lemma.  Then $G^* \cap
\calP_\kap * C(\kap^+) = G$ and for every
$\nu < \kap^+$ $i(p_\nu) \in G^*$.  Define
the elementary embedding $i^*: V[G] \to
M[G^*]$ by $i^* (\su a[G]) = (i (\su a))[G^*]$.  Then
$i^*|V = i$, $i^* (G)$ and, if $U^* = \{ A
\subset \kap| \kap \in i^* (A)\}$ then $U^*
= U(\kap, \alp, \emptyset)$.  So the
following diagram is commutative
$$\matrix{
&&M[G^*]&\cr
&i^*&&&\cr
&\nearrow&&\cr
V[G]&&\Bigg\uparrow&\ell\cr
&\searrow&&\cr
&&M^* \cong V[G]^\kap/U &\cr}$$
where $\ell([f]_{U^*}) = i^* (f) (\kap)$.

It remains to show that $\ell = id$.
Notice, that it is enough to prove that
every indiscernible in $\aC$ is of the form
$i^* (f) (\kap)$ for some $f \in V[G]$.
Examining the construction of $G^*$, it is
not hard to see that every indiscernible is
an element of the generic sequence $b_\del$
for some $\del \in k^{''} (N)$.

So indiscernibles are interpretations of
forcing terms with parameters in $k^{''}
(N)$.  But such elements can easily be
represented by functions on $\kap$ in
$V[G]$.  Let $\mu \in \aC$ and $\mu =
(t(\del) ) [G^*]$ for $\del = k(\del^*)$,
where $t$ is a term.  Let $\del^* =
[f]_{U(\kap, \alp)}$ for some $f: \kap \to
\kap$, $f \in V$.  Define a function $g \in
V[G]$ on $\kap$ as follows: $g(\tau) =
t(f(\tau))[G]$.  Then $i^* (g) (\kap) =
t(t(f) (\kap)) [G^*] = (t(\del))[G^*] =
\mu$.\hfill$\bigsquare$

Let us now turn to the construction of
$\ome_1$-preserving ideal.  Suppose that
$\alp^*$ is an $(\ome, \kap^+ + 1)$-repeat
point for $\aU$ in $V$.  Let $\alp^* \le
\alp < \alp^* + \kap^+$ be an ordinal.
Denote by $M^*$ the complete iteration of
$j_\kap^\alp (\aU)| ((j_\kap^\alp (\kap)$,
$(j_\kap^\alp(\alp^* + \kap^+))$ in
$N_\kap^\alp$.  Let $i^\alp$, $k^\alp$ be
the canonical elementary embeddings making
the diagram
$$\matrix{
&&N_\kap^\alp&\cr
&j_k^\alp & &\cr
&\nearrow &&\cr
V&&\Bigg\downarrow&k\cr
&i^\alp&&\cr
&\searrow&&\cr
&&M^\alp&\cr}$$
commutative.

As in the Main Lemma find in $V[G]$ a
$\calP(\kap, \alp$)-name of a
$M^\alp$-generic subset $\su G'$ of
$\calP_{i^*(\kap)} * C(i^\alp
(\kap^+))/\calP_{\kap + 1}$ so that each
$k^\alp (p_\nu^\alp) \in \su G'$ for $\nu <
\kap^+$, where $<p_\nu^\alp | \nu <
\kap^+>$ is the master condition sequence
for $U(\kap, \alp, \emptyset)$.  Let $\su
G^{''}$ be obtained from $\su G'$ by
removing all the information on the generic
subset of $C(i^\alp(\kap^+))$ except
$i^{\alp^{''}}(G \cap C(\kap^+))$.

We shall define a presaturated filter
$U^*(\kap, \alp)$ on $\kap$ in $V[G]$
extending $U(\kap, \alp)$ so that
\item{(i)} all generic ultrapowers of
$U^*(\kap, \alp)$ are generic extensions of
$M^\alp$;
\item{(ii)} every set $U(\kap, \alp,
\emptyset)$ is $U^*(\kap, \alp)$-positive.

Property (i) will insure that the
forcing for shooting clubs over
$i^\alp(\kap)$ will be $\kap$-closed
forcing.  Property (ii) is needed for
the iteration of forcings for shooting
clubs over $\kap$ in $V[G]$.

Denote $G \cap \calP_\kap$ by $\underline
G$ and $G \cap C(\kap^+)$ by $\oG$.  Let
$\oG (\nu)$ denote the $\nu$-th function
of $\oG$ for $\nu < \kap^+$.  Let $j^*$ be
the embedding of $V[G]$ into the ultrapower
of $V[G]$ by $U(\kap, \alp, \emptyset)$.
Let $< \alp_\nu|\nu < \kap^+ >$ be the list
of all the indiscernibles for
$i^\alp_\kap(\kap)$ of the complete
iteration used to define $M^\alp$.  For
every $\xi < \kap^+$ extend $\su G$ to $\su
G^*_\xi$ by adding to $\su G^{''}$
conditions $< \xi + \nu + 1, \kap,
\alp_\nu>$ for every $\nu < \kap^+$, i.e.
the $\xi + \nu + 1$ generic function moves
$\kap$ to $\alp_\nu$.

Define a filter $U^* (\kap, \alp)$ on
$\kap$ in $V[G]$ as follows:

$A \subseteq \kap$ belongs to $U^* (\kap,
\alp)$ iff for some $r \in G$ some $\su r'$
s.t. $\emptyset \vline \su r' \in \su
G^{''}$, for every $\su p \in \su C
(j_\kap^\alp (\kap^+))(0)$ s.t. in
$N_\kap^\alp[j^*(G)]$ $p \not\in j^*
(\oG(0))$, and for some $\su r^*$ s.t.
$\emptyset \vline \su r^* \in \su
G^*_{\xi(\su p)}$ in $M^\alp$
$$r \cup 1_{\calP (\kap, \alp)}\cup \su p
\cup \su r^*\buildrul {\calP_{i^\alp (\kap)
+ 1}} \under \vline \check \kap \in i^\alp
(\su A)$$
where $\xi(p) < \kap^+$ is the minimal s.t.
$p|\del \not\in j^*(G(0))$, for $\del$ the
$\xi(p)$ member of some fixed enumeration
of $j_\kap^\alp (\kap)$.

\proclaim Claim 1.  $U(\kap, \alp,
\emptyset) \supseteq U^* (\kap, \alp)$.

\pr Let $A \in U (\kap, \alp, \emptyset)$,
then for some $\nu < \kap^+$ some $r$, $<
\emptyset, T > \in \calP(\kap, \alp)$, in
$N_\kap^\alp$
$$r \cup \{< \emptyset, T > \} \cup \su
p_\nu \vline \check \kap \in j_\kap^\alp
(\su A)\ .$$
Let us split $\su p_\nu$ into three parts
$\su p_1 = \su p_\nu \mid
\calP_{j_\kap^\alp}$, $\su p_2 = \su p_\nu
\cap C(j_\kap^\alp (\kap^+))(0)$ and $\su
p_3$ = the rest of $\su p_\nu$.  Then,
using $k^\alp$, in $M^\alp$
$$r \cup \{ < \emptyset, T > \} \cup k^\alp
(\su p_\nu) \vline \check \kap \in i^\alp
(\su A)\ .$$
Extend $\su p_2$ to some $\su p'_2$, still
remaining in the master condition
sequence, in order to make its domain
above all the indices of the generic
sequences listed in $\su p_3$.  Extend
$p'_2$ to some $p^{''}_2 \in C(j_\kap^\alp
(\kap^+))(0)$ which is incompatible with a
member of the master condition sequence.
Then, using $k^\alp$, in $M^*$ 
$$r \cup \{ < \emptyset, T > \} \cup k^\alp
(\su p_1) \cup k^\alp (p^{''}_2) \cup
k^\alp (\su p_3) \vline \check \kap \in
i^\alp (\su A)\ .$$
By the definition of $U^* (\kap, \alp)$, it
implies that $A$ is a $U^* (\kap,
\alp)$-positive set.

So every member of $U(\kap, \alp,
\emptyset)$ is $U^* (\kap, \alp)$-positive.
The fact that $U(\kap, \alp, \emptyset)$ is
an ultrafilter completes the proof of the
claim.

\proclaim Claim 2.  $U^*(\kap, \alp)$ is
normal precipitous cardinal preserving
filter and its generic ultrapower is a
generic extension of $M^\alp$.

\pr Let $U = U^*(\kap, \alp) \cap
V[\underline G, \overline G(0)]$.  The
arguments of Theorem 1  show that a
generic ultrapower by $U$ is isomorphic to
a generic extension of the complete
iteration of $N_\kap^\alp$ with
$j_\kap^\alp(\vec U) |
j_\kap^\alp(\kap)$ above $\kap$.  Actually
the forcing with $U$ is isomorphic to
$\calP(\kap, \alp)$ followed by
$C(j_\kap^\alp (\kap^+))(0)$ (in the sense of
this iteration) over $V[\underline G,
\oG(0)]$.  Clearly, the generic
function form $j_\kap^\alp(\kap^+)$ into
$j_\kap^\alp(\kap^+)$  produced by this
forcing is incompatible with $j^*(\oG(0))$
which belongs to $V[\underline G, \oG(0)]$.

Now the forcing with $U^* (\kap, \alp)$
over $V[G]$ does the following: First it
picks $\su p \in \su C
(j_\kap^\alp(\kap^+)) (0)$ incompatible with
$j^* (\oG(0))$ and then $G^*_{\xi(\su p)}$
is added.  $G^* _{\xi (\su p)}$ insures
that a generic ultrapower is a generic
extension of $M^\alp$.  So the forcing with
$U^* (\kap, \alp)$ is isomorphic to
$\calP(\kap, \alp)$ followed by a portion
 of $C(i^\alp (\kap^+))$, which is
$\kap$-closed.\hfill$\bigsquare$

Force over $V[G]$ with the forcing $Q_\kap$
which is the Backward-Easton iteration of
the forcings adding $\del^+$-Cohen subsets
to every regular $\del < \kap$ with
$0^{\vec U} (\del) > 0$.
Fix a generic subset $H$ of $Q_\kap$.  All
the filters
$U^* (\kap, \alp)$ extend in the obvious
fashion in $V[G,H]$.  Let us use the same
notations for the extended filters.

Set $F = \cap \{U^* (\kap, \alp)|\alp^* \le
\alp \le \alp^* + \kap^+\}$.  Then $F$ is a
$\le \kap$-preserving filter in $V[G, H]$
and forcing with it is isomorphic to $\calP
(\kap, \alp)$ (for some $\alp, \alp^* \le
\alp \le \alp^* + \kap^+)$ followed by
$\kap$-closed forcing.  Let us shoot clubs
through  elements of $F$, then through the
sets  of generic points and so on, as 
was done in [G 3,4].  Denote this forcing
by $B$.  Let $R$ be its generic subset over
$V[G, H]$.  We shall show that $NS_\kap$
is $\ome_1$-preserving ideal in $V[G,H,R]$.
The forcing with $NS_\kap$ consists of
two parts:
\item{(a)} embedding of $B$ into
$\calP(\kap, \alp)$;
\item{(b)} $(\kap$-closed forcing) $*
(i^\alp (B) / i^{'' \alp} (R))$.

By the choice of $i^\alp (\kap)$, the
forcing $i^\alp (B)$ is a shooting club
through sets containing a club (the club of
indiscernibles for $i^\alp(\kap))$.  So it
is $\kap$-closed forcing and part (b)
does not cause any problem.

Let us examine part (a) and show that
this forcing preserves $\ome_1$.  Recall
that $B$ is the direct limit of $< B_\bet |
\bet < \kap^+>$ where each $B_\bet$ is of
cardinality $\kap$ and for a limit $\bet$,
$B_\bet = $ the direct limit of
$B_{\bet'}(\bet' < \bet)$, if $cf \bet =
\kap$ and $B_\bet$ = the inverse limit of
$B_{\bet'} (\bet' < \bet)$ otherwise.  The
forcing of (a) is
$$\calP = \{ \pi \in V[G,H]| \ \hbox{for
some}\ \bet < \kap^+\enskip \pi \ \hbox{ is an
embedding of}\ B_\bet\ \hbox{into}\
\calP(\kap, \alp)\}\ .$$
For $\pi_1, \pi_2 \in \calP$ let $\pi_1 \ge
\pi_2$ if $\pi_1| \dom \pi_2 = \pi_2$.

\proclaim Claim.  $NS_\kap$ is
an $\ome_1$-preserving ideal in $V[G,H,R]$.

\pr Suppose otherwise.  Then for some
$\alp^* \le \alp \le \alp^* + \kap^+$ some
condition in the forcing $\calP(\kap,
\alp)^*$ (the forcing isomorphic to the
adding of $\kap^+$-Cohen subsets of
$\kap^+$) $* \calP$ over $V[G,H]$ forces
$\ome_1$ to collapse.  Let us assume
that the empty condition already forces
this.  Consider the case when
$\vline_{\calP(\kap, \alp)} cf \check \kap
= \check \ome$.  The remaining cases are
simpler.

Let $S$ be a $V[G,H]$-generic subset of
$\calP(\kap, \alp) * (\kap^+$-Cohen subsets
of $\kap^+)$.  Denote $V[G,H,S]$ by $\oV$.
Let $\su f$ be a $\calP$-name in $\oV$ of a
function from $\ome$ to $\ome_1^{\oV}$.
Pick an elementary submodel $N$ of $<
H(\lam), \eps, B, \calP, \kap, \su f > $,
for $\lam$ big enough, satisfying the
following three conditions:
\item{(1)} $|N| = \kap$;
\item{(2)} $N \supseteq H^\oV (\kap)$;
\item{(3)} $N \cap \kap^+ = \del$ for some
$\del$ s.t. $cf^{V[G, H]}\del = \kap$.

Then $B \cap N = B_\del$.  Also $B_\del$ is
a direct limit of $< B_{\del'} | \del' <
\del>$.  Pick in $V[G,H]$ a cofinal
sequence $< \del_\bet|\bet < \kap>$ to
$\del$ and in $\oV$ a cofinal sequence $<
\tau_n|n<\ome>$ to $\kap$.  Consider the
subsets of $\kap$, $< \su A_\bet | \bet <
\kap >$ s.t. $B_{\del_\bet + 1} /
B_{\del_\bet}$ is the forcing for shooting
club into $\su A_\bet$.  Assume for
simplicity that all $A_\bet$'s are in $V$.
Let $A = \Del_{\bet < \kap} A_\bet = \{
\gam< \kap|$ for every $\bet < \gam, \quad
\gam \in A_\bet\}$.  The $A_\bet$'s and $A$
are in $\cap \{ U(\kap, \gam)| \alp^* \le
\gam \le \alp^* + \kap^+\}$.  Pick $\pi_0
\in \calP \cap N$ deciding the value of
$\su f(0)$.  Without loss of generality $\dom
\pi_0 = B_{\del_{\bet_0} + 1}$ for some
$\bet_0 < \del$.  Let $n_0$ be the least $n
< \ome$ such that $\del_{\tau_n} > \bet_0$.
Denote $\del_{\tau_{n_0}}$ by $\eps_0$.  Let
$C_0$ be the generic club through
$A_{\bet_0}$ defined by $\pi_0$ and $S$.
As in Lemma 3.6 [G4], it is possible to
find an element of the generic sequence to
$\kap$, $\tau(0) \in A - \tau_{n_0}$ such
that $C \cap \tau$ is a $V[G| \tau(0)$, $H|
\tau(0)]$-generic club through $A \cap
\tau$.  Choosing $\tau(0)$ 
more carefully, it is possible to also satisfy
the following $A \cap \tau(0) = \Del_{\bet
< \tau(0)} (A_\bet \cap \tau(0))$.  Now pick
$\pi'_0 \in \calP \cap N$ to be an
extension of $\pi_0$ with domain
$B_{\del_r}$ such that the clubs of $\pi'_0$
intersected with $\tau(0)$ are $V[G| \tau
(0)$, $H| \tau(0)]$-generic for $B_\del|
\tau(0)$.  It is possible since $N$
satisfies condition (2).

Now find in $N$ an extension $\pi_1$ of
$\pi'_0$ deciding $\su f(1)$.  Define
$\pi'_1$ and $\tau(1)$ as above.  Continue
the process for all $n < \ome$.  Finally
set $\pi = \bigcup\limits_{n < \ome}
\pi_n$.  It is enough to show that $\pi \in
\calP$, since then $\pi \vline \su f \in
V_1$.  Let us prove that $\pi$ and $S$
produce a $V[G, H]$-generic subset of
$B_\del$.  Suppose that $\calD \in V[G,H]$
is a dense subset of $B_\del$.  Then $X =
\{ \tau < \kap| \calD \cap H(\tau)$ is a
dense subset of $B_{\del_r} | \tau\}$
contains a club in $V[G,H]$.  Since the
generic sequence to $\kap$ is almost
contained in every club of $V[G,H]$, for
some $n< \ome$, $\tau(n) \in X$.  But then
the generic subset produced by $\pi'_n$
intersects $\calD$.  So the same is
true for $\pi$.\hfill$\bigsquare$

Let us now turn to successor cardinals.  We
would like to make $NS_\kap$ $\ome_1$
preserving for $\kap = \mu^+$ for a
regular $\mu > \aleph_1$.  It is possible
to use the model constructed above,
collapse $\kap$ to $\mu^+$ and apply
the results of Section 1.  But an $(\ome, \kap^+ +
1)$-repeat point was used in the
construction of the model.  It turns out
that an $(\ome, \mu + 1)$-repeat point
suffices for $NS_{\mu^+}$ and an $(\ome,
\mu)$-repeat point for $NS_{\mu^+} | \{\alp
< \mu^+ | cf \alp < \mu\}$.  On the other
hand, precipitousness of
$NS_{\mu^+}^{\aleph_0}$ implies an $(\ome,
\mu)$-repeat point, by [G4].

Let us preserve the notations used above.
Assume that $\mu < \kap$ is a regular
cardinal and some $\alp^* < 0^{\vec
 U} (\kap)$ is an $(\ome, \mu +
1)$-repeat point.  Let $G$ be a generic
subset of $\calP_\kap * C(\kap^+)$.  Over
V[G] instead of the forcing $Q_\kap$, in
the previous construction, use $\Col (\mu,
\kap)$ the Levy collapse of all the
cardinals $\tau, \mu < \tau < \kap$ on
$\mu$.  Let $H$ be a generic subset of
$\Col(\mu, \kap)$.  Denote by $H(\tau)$ the
generic function from $\mu$ on $\tau$ where
$\tau \in (\mu, \kap)$.

Now the forcing for shooting clubs should
come.  In order to prevent collapsing
cardinals by this forcing, $j(\kap)$ was
made a limit of $\kap^+$ indiscernibles for
the measures $\{U(\kap, \alp)| \alp^* \le
\alp < \alp + \kap^+\}$.  But now we have
only $\mu$ measures.  So the best we can do
is to make $j(\kap)$ a limit of $\mu$
indiscernibles and then its cofinality in
$V$ will be $\mu < \kap$.  It looks 
slightly paradoxical since usually $cf
j(\kap) = \kap^+$, but it is possible by
[G5].  In order to explain the idea of [G5]
which will be used here let us give an
example of a precipitous ideal $I$ on
$\ome_1$ so that:
$$\buildrul I^+ \under \vline ''cf^V
(j(\check \ome_1)) = \check \ome''\
.\leqno(*)$$

\subheading{Example}  Suppose that $\kap$
is a measurable cardinal.  Let $U$ be a
normal measure on $\kap$ and $j: V \to N
\cong V^\kap/U$ the canonical embedding.
Let $\calP_\kap$ be the Backward-Easton
iteration of the forcings $C(\alp^+)$ for
all regular $\alp < \kap$.  Let $G_\kap *
H_\kap$ be a generic subset of $V$.
Collapse $\kap$ to $\ome_1$ by the Levy
collapse.  Let $R_\kap$ be $V[G_\kap *
H_\kap]$-generic subset of Col$(\ome_1,
\kap)$.  Denote $V[G_\kap * H_\kap *
R_\kap]$ by $V_1$.  We shall define a
precipitous ideal satisfying $(*)$ in $V_1$.

Let $j_0 = j$, $N_0 = N$, $\kap_0 = \kap$,
$U_0 = U$.  Set $N_1 \cong N_0^{j_0(\kap)}
/ j_0 (U)$, $\kap_1 = j_0 (\kap)$, $U_1 =
j_0(U)$ and let $j_1: N_0 \to N_1$ be the
canonical embedding.  Continue the
definition 
for all $n < \ome$.  Set $j_\ome,
N_\ome$ to be the direct limit of $<
j_n, N_n|n < \ome >$.  Then $j_\ome(\kap) =
\bigcup\limits_{n < \ome} \kap_n$.  Set
$\kap_\ome = j_\ome(\kap)$.  Notice, that
$\cup (j^{''} (\kap^+)) = (\kap_1^+)^{N_0}$
and $\cup(j_{n+1} ((\kap_n^+)^{N_n})) =
(\kap_{n+1}^+)^{N_n + 1}$.  So
$\cup(j_\ome^{''}(\kap^+)) =
(\kap_\ome^+)^{N_\ome}$.

Define a filter $U^*$ in $V[G_\kap * H_\kap
* R]$ as follows:
\hb
$A \in U^*$ iff for some $r \in G_\kap *
H_\kap * R$ for some $n < \ome$ in $N_\ome$
$$r \cup p_n \vline \check \kap \in j_\ome
(\su A)$$
where $p_n$ is the name of condition in the
forcing $C(j_\ome(\kap^+))$ defined as
follows:\hb
let $\su h$ be the name of the generic
function from $\ome$ onto $(\kap^+)^V$ in
Col$(\ome, j_\ome(\kap))/G_\kap * H_\kap *
R_\kap$.  Set $p_n = \{ < j_\ome (\su h
(0)), \kap, \kap_1>$, $< j_\ome (h(\su 1)),
\kap, \kap_2> \nek < j_\ome (\su h (n)),
\kap, \kap_{n+1} > \}$.  The meaning of the
above is that the value on $\kap$ of
$j_\ome(h(m))$-th function from
$j_\ome(\kap)$ to $j_\ome(\kap)$ is forced
to be $\kap_{m+1}$.

It is not difficult to see now that $U^*$ is
a normal precipitous ideal on $\ome_1$ and a
generic ultrapower with it is isomorphic to
a $V[G_\kap ^* H_\kap * R_\kap]$-generic
extension of $N_\ome[G_\kap, H_\kap,
R_\kap]$.  Also for a generic embedding
$j^*$, $j^*| V = j_\ome$.  Hence $j^*
(\kap) = \kap_\ome$ which is of cofinality
$\ome$ in $V$.

As in [J-M-Mi-P], it is possible to extend
$U^*$ to the closed unbounded filter with
the same property.  Using the Namba
forcing, it is possible to construct a
precipitous ideal satisfying $(*)$ on
$\aleph_2$.  Starting from a measurable
which is a limit of measurable, it is
possible to build such an ideal over an
inaccessible or even measurable.  Since
then it is possible to change the
cofinality of the ordinal of cofinality
$\kap^+$ to $\ome$ in $N$ and that is what
was needed  to catch all $\kap_n$'s
in the above construction. We do not know
if one measurable is sufficient for a
precipitous ideal satisfying $(*)$ over
$\kap > \aleph_2$.

Note also that by Proposition 1.5 if $I$ is
$\lam$-preserving, then $cf^V j(\kap) \ge
\lam$.  In particular, if $I$ is
presaturated then $cf^V j(\kap) = \kap^+$.

Let us now return to the construction of
$NS_{\mu^+}$\ $\ome_1$-preserving.  Let
$\alp^* \le \alp < \alp^* + \mu$ be an
ordinal.  Define $M^*$ as in the
construction of $NS_\kap$\ 
$\ome_1$-preserving for an inaccessible
$\kap$.  Now $cf^V i^\alp (\kap)$ will be
$\mu$.  We shall define, in $V[G,H]$, $U^*
(\kap, \alp)$ extending $U(\kap, \alp)$
satisfying conditions (i) and (ii)
from page 23.  To do so simply
combine the definition there with the
definition of the example.  We leave the
details to the reader.  The rest of the
construction does not differ from an
inaccessible cardinal case.

The above results give  equiconsistency for
$NS_\kap|$ (singular) or something close to
equiconsistency for $NS_\kap$, but for
$\kap > \aleph_2$.  For $\kap = \aleph_2$,
we do not know if the assumption of the
existence of an $(\ome, \ome_1 + 1)$-repeat
point (or of an $(\ome, \ome_1)$-repeat point
for $NS_{\aleph_2}^{\aleph_0}$) can be
weakened.  By [G3], a measurable is
sufficient for the precipitousness of
$NS_{\aleph_2}^{\aleph_0}$ and a measurable
of order 2 for $NS_{\aleph_2}$.  Let us
show that $\ome_1$-preservingness requires
stronger assumptions.

\proclaim Lemma 2.5.  Suppose that $\kap >
\aleph_1$ is a regular cardinal, $2^\kap =
\kap^+$, $2^{\aleph_0} = \aleph_1$.  $I$ is
a normal $\ome_1$-preserving ideal over
$\kap$ so that $\{\alp < \kap|cf \alp =
\ome\} \not\in I$.  Then $\exists \tau
0(\tau) \ge 2$ in the core model.

\pr Without loss of generality, assume that
$\neg \exists \alp 0(\alp) = \alp^{++}$.
Denote by $\calK(\vec
\calF)$ the core model with the maximal
sequence $\vec F$.  We refer
to Mitchell papers [Mi 1,2] for the
definitions and properties of
$\calK(\vec \calF)$ that we
are going to use.

The set $A = \{ \alp < \kap | \alp$ is
regular in $\calK(\vec
\calF)$ and of cofinality $\ome$ in $V\}$
is $I$-positive.  If the set $A^* = \{ \alp
\in A| \alp$ is not measurable in
$\calK(\vec \calF)\}$ is
bounded in $\kap$, then $0(\kap) \ge 2$.
Suppose otherwise.  Let $j:V \to M$ be a
generic elementary embedding so that the
set $\{\alp < \kap|cf \alp = \ome\}$
belongs to a generic ultrafilter $G_I$.
Then $j(A^*)$ is unbounded in $j(\kap)$.
Pick the mimimal $\alp \in j(A^*) - \kap$.

\proclaim Claim.  $cf^{\calK(\vec
 \calF)} (\alp) =
(\kap^+)^{\calK(\vec
\calF)}$.

\pr Suppose otherwise.  Then, since $\alp$
is regular in $\calK(\aF)$ and $j|
\calK(\aF)$ is an iterated ultrapower of
$\calK(\aF)$ by $\aF$, $\alp$ is a limit
indiscernible of this iteration.  By
Proposition 1.5.  $^\ome M \cap V[G_I]
\subseteq M$.  Since $\alp$ is not a
measurable in $\calK(j(\aF))$, it implies
that $cf ^{\calK(\aF )} \alp > \ome$.  But
then, using $^\ome M \cap V[G_I] \subseteq
M$ and the arguments of [Mi 2] and [G4], we
obtain some $\tau$ with $0^{\aF} (\tau) \ge
(\ome_1)^{ \calK(\aF)}$. $\square$ of the
claim.

By [Mi 2], $(\kap^+)^{\calK(\aF)} =
(\kap^+)^V$.  But in $V[G_I]$, $cf \alp =
\ome$ and hence $cf (\kap^+)^V = \ome$
which contradicts Proposition
1.5.\hfill$\bigsquare$

\proclaim Lemma 2.6.  Suppose that $\kap >
\aleph_1$ is a regular cardinal, $2^\kap =
\kap^+$, $2^{\aleph_0} = \aleph_1$ and $I$
is a normal $\ome_1$-preserving ideal.  If
there exists an $\ome$-club $C$ so that
every $\alp \in C$ is regular in
$\calK(\aF)$, then $\exists \tau 0^{\aF}
(\tau) \ge \ome_1$ in $\calK(\aF)$.

The proof is similar to Lemma 2.5; just
consider the $\ome_1$-th member of $j(C) -
\kap$.

\proclaim Theorem 2.7.  Assume $GCH$, if
$NS_{\aleph_2}^{\aleph_0}$ is
$\ome_1$-preserving then $\exists \tau$
$0^{\aF} (\tau) \ge \ome_1$ in $\calK(\aF)$.
The proof follows from Lemma 2.6.
\bigskip

\references{55}

\ref{[A]} U. Avraham, Isomorphism of
Aronszajn trees, Ph.D. Thesis, Jerusalem,
1979.
\smallskip
\ref{[B]} J. Baumgartner, Independence
results in set theory, Notices Am. Math.
Soc. 25 (1978), A 248-249.
\smallskip
\ref{[B-T]} J. Baumgartner and A. Taylor,
Ideals in generic extensions.  II, Trans.
of Am. Math. Soc. 271 ( 1982), 587-609.
\smallskip
\ref{[G-J-M]} F. Galvin, T. Jech and M.
Magidor, An ideal game, J. of Symbolic
Logic 43 (1978), 284-292.
\smallskip
\ref{[G1]} M. Gitik, Changing cofinalities
and the nonstationary ideal, Israel J. of
Math. 56 (1986), 280-314.
\smallskip
\ref{[G2]} M. Gitik, $\neg SCH$ from
$0(\kap) = \kap^{++}$, Ann. of Pure 43
(1989), 209-234.
\smallskip
\ref{[G3]} M. Gitik, The nonstationary
ideal on $\aleph_2$, Israel J. of Math.
48 (1984), 257-288.
\smallskip
\ref{[G4]} M. Gitik, Some results on the
nonstationary ideal, Israel Journal of
Math.
\smallskip
\ref{[G5]} M. Gitik, On generic elementary
embeddings, J. Sym. Logic 54(3) (1989),
700-707.

\smallskip
\ref{[H-S]}  L. Harrington and S. Shelah,
Equiconsistency results in set theory,
Notre Dame J. of Formal Logic 26(2) (1985), 
 178-188.
\smallskip
\ref{[J-M-Mi-P]} T. Jech, M. Magidor, W.
Mitchell and K. Prikry, Precipitous ideals,
J. Sym. Logic 45 (1980), 1-8.
\smallskip
\ref{[Mi 1]} W. Mitchell, The core model
for sequences of measures I, Math. Proc.
Camb. Phil Soc. 95 (1984), 229-260.
\smallskip
\ref{[Mi 2]} W. Mitchell, The core model
for sequence of measures II, to appear.

\smallskip
\ref{[Mi 3]} W. Mitchell, Indiscernibles,
skies and ideals, Contemp.  Math. 31
(1984), 161-182.
\smallskip
\ref{[S1]} S. Shelah, Proper forcing,
Springer-Verlag Lecture Notes in Math,
940 (1982).
\smallskip
\ref{[S2]} S. Shelah, Some notes on
iterated forcing with $2^{\aleph_0} >
\aleph_2$, Notre Dame J. of Formal Logic.

\end